\documentclass[a4paper,12pt]{article}
\usepackage{amssymb,amsmath}
\usepackage[all]{xy}
\usepackage{mathabx}
\usepackage{empheq}
\newtheorem{theorem}{Theorem}[section]
\newtheorem{lemma}[theorem]{Lemma}
\newtheorem{proposition}[theorem]{Proposition}
\newtheorem{corollary}[theorem]{Corollary}
\newcommand{\udots}{\mathinner{\mskip1mu\raise1pt\vbox{\kern7pt\hbox{.}}
   \mskip2mu\raise4pt\hbox{.}\mskip2mu\raise7pt\hbox{.}\mskip1mu}}

\newenvironment{definition}[1][Definition]{\begin{trivlist}
\item[\hskip \labelsep {\bfseries #1}]}{\end{trivlist}}

\newenvironment{def-prop}[1][Definition-Proposition]{\begin{trivlist}
\item[\hskip \labelsep {\bfseries #1}]}{\end{trivlist}}
\newenvironment{remark}[1][Remark]{\begin{trivlist}
\item[\hskip \labelsep {\bfseries #1}]}{\end{trivlist}}

\newcommand{\nce}{\nu\{\# m\mathrm{-bad}\}}
\newcommand{\ce}{\{\# m-\mathrm{bad}\}}
\newcommand{\nun}{\nu\{\#j-\mathrm{bad}\}}
\newcommand{\ndo}{\nu\{\mathrm{positions\ of\ }j'\mathrm{s} \}}
\newcommand{\ntr}{\nu\{\#f\}}
\newcommand{\ncu}{\nu\{f \mathrm{\ to\ the\ right}\}}
\newcommand{\nci}{\nu\{\mathrm{depth\ }m\mathrm{\ and\ }j \}}
\newcommand{\nse}{\nu\{m \mathrm{\ far\ from\ bottom}\}}
\newcommand{\nsi}{\nu\{\mathrm{ min\ }m-\mathrm{bad}\} }
\newcommand{\noc}{\nu\{\#m, j \mathrm{\ after\ min.\ }m-\mathrm{bad}\}}
\newcommand{\nnu}{\nu\{\mathrm{function\ of\ } m-\mathrm{bads}\} }
\newcommand{\nue}{\{\mathrm{function\ of\ } m-\mathrm{bads}\} }
\newcommand{\ndi}{\nu\{\mathrm{max.\ } j-\mathrm{bad} \}}
\newcommand{\non}{\nu\{m,j  \mathrm{\ equal\ to\ left}\}}

\newcommand{\CE}{{\cal E}}

\title{Presentation of right-angled Soergel categories by generators and relations}
\author{Nicolas Libedinsky}

\begin{document}

\maketitle

\begin{abstract}
Soergel bimodule category \textbf{B} is a categorification of the Hecke algebra of a Coxeter system $(W,\mathcal{S})$. We find a presentation of \textbf{B} (as a  tensor category) by generators and relations when $W$ is a right-angled Coxeter group.
\end{abstract}

\section{Introduction}

In 1979 \cite{KL1}, Kazhdan and Lusztig defined the Kazhdan-Lusztig polynomials, and this gave rise to what now is  known as the Kazhdan-Lusztig theory. They stated two major conjectures involving these polynomials in \cite{KL1} and \cite{KL2} : Kazhdan-Lusztig conjecture (in representation theory of Lie algebras) and Kazhdan-Lusztig positivity conjecture (in algebraic combinatorics).  The first conjecture was proved for Weyl groups by Beilinson and Bernstein in \cite{BB},  by Brylinski and Kashiwara in \cite{BK} and later by Soergel \cite{S0}. The second conjecture has been proved for Weyl or affine Weyl groups in \cite{KL2} and in some other cases by Haddad \cite{H} and Dyer \cite{Dy}.

In 1980 \cite{Lu}, Lusztig stated a central conjecture in representation theory, known as Lusztig conjecture. This conjecture is about the characters values of the irreducible representations of algebraic groups in positive characteristic. This conjecture is solved only for large characteristics. 

Let us consider $(W,\mathcal{S})$ a Coxeter system and $\mathcal{H}$ its Hecke algebra. In 1992 \cite{S1}, Soergel categorified $\mathcal{H}$. This means that he defined a tensor category $\mathbf{B}$ (that depends on a field $k$ and on a representation of $W$) and an isomorphism of rings $\CE$
 from $\mathcal{H}$ to the split Grothendieck group of
  $\mathbf{B}$. He has then stated a conjecture that links, via $\CE$, the  Kazhdan-Lusztig basis elements in $\mathcal{H}$ with the indecomposable elements of  $\mathbf{B}$. This conjecture implies  Kazhdan-Lusztig positivity conjecture, and when the characteristic of $k$ is larger that the coxeter number $h$ of $W$, it implies a part of Lusztig conjecture.

Soergel introduced this category in the case of Weyl groups to make a link between the BGG  category $\mathcal{O}$ of a semisimple complex Lie algebra and the perverse sheaves over the grassmanians. This link with geometry allowed him to prove his conjecture in the Weyl group case. He deduced in \cite{S0} the proof of Kazhdan-Lusztig conjecture we have mentioned. 

The following diagram is a summery of the implications 
\vspace{0.4cm}

$
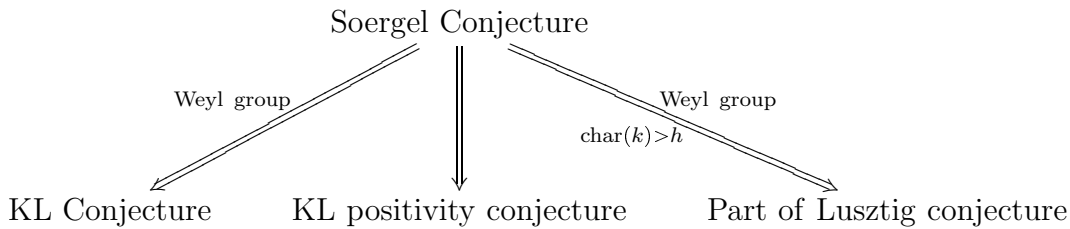
\xymatrix{
  \vspace{-5cm} &  \mathrm{Soergel\ Conjecture}   \ar@{=>}[ldd]_{\mathrm{Weyl\ group\ \ }} \ar@{=>}[dd] \ar@{=>}[rdd]^{\ \ \mathrm{Weyl\ group}}_{\mathrm{char}(k)>h\ \ }&  \\
&&\\
   \mathrm{KL\ Conjecture }  & \mathrm{KL\ positivity \ conjecture  } & \mathrm{Part\ of\ Lusztig \ conjecture }
  }
$

\vspace{0.4cm}
In this paper, we are concerned with the case where $(W,\mathcal{S})$ is a right-angled Coxeter system. This means that $m(s,r)=2$ or $\infty$ for all $s,r\in\mathcal{S}.$ In this case, we find a presentation as tensor category of $\mathbf{B}$ (called right-angled Soergel category) by generators and relations. Our main motivation to do this is the following. If we can extend this result to any Coxeter system we expect to be able to do the following  : to reprove Kazhdan-Lusztig conjecture in the spirit of \cite{S0} but (from our perspective) in a more natural way and  to reprove a part of Lusztig conjecture  (comparison between quantum group and algebraic group) in an essentially different way that the proof in \cite{AJS}.

To extend the result of this paper  to the symmetric group would be a step forward in the program of calculating Khovanov-Rozansky link homology (a categorification of the HOMFLYPT polynomial) defined in \cite{K}  via Soergel category \textbf{B}. Some results in this direction are given in the papers \cite{Ra} and \cite{WW}.

We should say that the right-angled case is very rich in topological and geometrical terms. For example, in the introduction of \cite{D} Davis remarks that the right-angled case is sufficient for the construction of most examples of interest in geometric group theory.  

This paper is divided as follows: in sections \ref{section2} and \ref{section3} we have the basic definitions and preliminaries. From section \ref{section4} to section \ref{section8} we prove the central theorem \ref{bacanz}. At the beginning of section \ref{section4} we make a summary of the proof.

\section{Notations and preliminaries}\label{section2}

\subsection{Basic definitions}

Let $(W,\mathcal{S})$ be a not necessarily finite Coxeter system (with $\mathcal{S}$ a finite set) and $\mathcal{T}\subset W$ the set of reflections in $W$, \textit{i.e.} the orbit of $\mathcal{S}$ under conjugation. Let $k$ be an infinite field of characteristic different from $2$ and $V$ a finite dimensional $k$-representation of $W$. For $w\in W$, we denote by $V^w\subset V$ the set of $w-$fixed points. The following definition might be found in \cite{S3} :

\vspace{0.5cm}  \begin{definition}
 By a \textit{reflection faithful representation}  of $(W,\mathcal{S})$ we mean a faithful, finite dimensional representation $V$ of $W$ such that, for each $w\in W$, the subspace $V^w$ is an hyperplane of $V$ if and only if $w\in\mathcal{T}.$
\end{definition}

 From now on, we consider $V$ a reflection faithful representation of $W$. If $k=\mathbb{R}$, by the results of  \cite{lib}, all the results in this paper will stay true if  we consider $V$ to be the geometric representation of $W$ (still if this representation in not always reflection faithful).

Let $\widehat{R}=R(V)$\label{d1} be the algebra of regular functions on $V$. The 
action of $W$ on $V$ induces an action on $\widehat{R}$. For $s\in\mathcal{S}$ consider the $(\widehat{R},\widehat{R})-$bimodule $\widehat{\theta}_s=\widehat{R}\otimes_{\widehat{R}^s} \widehat{R}$, where $\widehat{R}^s$ is the subspace of $\widehat{R}$ stabilized by $s$.

\vspace{0.5cm}  \begin{definition}
 Soergel's category $\mathbf{B}(W,V)=\mathbf{B}$ is the full subcategory of all $(\widehat{R},\widehat{R})-$bimodules with objects the finite direct sums of direct summands of bimodules of the type $\widehat{\theta}_{s_1}\otimes_{\widehat{R}}\widehat{\theta}_{s_2}\otimes_{\widehat{R}}\cdots\otimes_{\widehat{R}} \widehat{\theta}_{s_n}$ for $(s_1,\ldots, s_n)\in \mathcal{S}^n.$ 

For simplicity, we will denote by $\widehat{\theta}_{s_1}\widehat{\theta}_{s_2}\cdots  \widehat{\theta}_{s_n}$ the $(\widehat{R},\widehat{R})-$bimodule $\widehat{\theta}_{s_1}\otimes_{\widehat{R}}\widehat{\theta}_{s_2}\otimes_{\widehat{R}}\cdots\otimes_{\widehat{R}} \widehat{\theta}_{s_n}\cong \widehat{R}\otimes_{\widehat{R}^{s_1}}\widehat{R}\otimes_{\widehat{R}^{s_2}}\cdots \otimes_{\widehat{R}^{s_n}}\widehat{R}$
\end{definition}

More details about this category can be found in the papers by Dyer \cite{Dy2}, \cite{Dy3}, by Fiebig \cite{Fie1}, \cite{Fie2}, \cite{Fie3}, by myself \cite{L}, \cite{lib}, by Soergel \cite{S0}, \cite{S1}, \cite{S2}, \cite{S3} and by Williamson \cite{Wi}, being \cite{S3} the central and more complete reference from our point of view.

We will use an auxiliary subcategory of $\mathbf{B}$, where we do not consider the direct summands :

\vspace{0.5cm}  \begin{definition}\label{d155}
$ \widehat{\mathbf{B}}$ is the full subcategory of $\mathbf{B}$ with objects the finite direct sums  of bimodules of the type $\widehat{\theta}_{s_1}\widehat{\theta}_{s_2}\cdots \widehat{\theta}_{s_n}$ for $(s_1,\ldots, s_n)\in \mathcal{S}^n.$
\end{definition}

\subsection{Stability of $V'$}\label{stab}

Let $\widehat{x}_s\in V^*$ be an equation of the hyperplane $H_s$ fixed by $s\in\mathcal{S}$. We have a decomposition $\widehat{R}\simeq \widehat{R}^s\oplus \widehat{x}_s\widehat{R}^s,$ corresponding to $$\widehat{R}\ni p=\frac{p+s\cdot p}{2}+\frac{p-s\cdot p}{2}.$$ We define $\widehat{P}_s(p)=(p+s\cdot p)/2$, $\widehat{I}_s(p)=(p-s\cdot p)/2$\label{d2} and $\widehat{\partial}_s(p)=(p-s\cdot p)/2\widehat{x}_s.$

As $W$ acts trivially over $$(\sum_{s\in \mathcal{S}}k\widehat{x}_s)^{\perp} = \bigcap_{s\in \mathcal{S}}H_s,$$
we have that $V'=\sum_{s\in \mathcal{S}}k\widehat{x}_s$ is stabilized by $W$, so we deduce the important fact that  $\widehat{P}_t(\widehat{x}_s), \widehat{I}_t(\widehat{x}_s)\in V'. $

\subsection{Some morphisms}
We define four morphisms in $\widehat{\mathbf{B}}$ :

\begin{displaymath}\label{d3}
\begin{array}{lll}
     \hspace{-4.7cm}   \hspace{3.2cm} \widehat{j}_r : \widehat{\theta}_r\widehat{\theta}_r &\rightarrow &  \widehat{\theta}_r \\
\hspace{-4.7cm}\widehat{R}\otimes_{\widehat{R}^r}\widehat{R}\otimes_{\widehat{R}^r}\widehat{R}\ni p_1\otimes p_2\otimes p_3 &\mapsto & p_1\widehat{\partial}_r(p_2)\otimes p_3 ;
\end{array}
\end{displaymath}

\begin{displaymath}
\begin{array}{lll}
\hspace{-4.3cm}        \hspace{1.8cm} \widehat{m}_r :\widehat{\theta}_r &\rightarrow &  \widehat{R} \\
\hspace{-4.3cm}\widehat{R}\otimes_{\widehat{R}^r}\widehat{R}\ni p_1\otimes p_2 &\mapsto & p_1p_2 ;
\end{array}
\end{displaymath}

\vspace{0.3cm}

\begin{displaymath}
\begin{array}{lll}
      \hspace{1.2cm}  \widehat{\alpha}_r :\widehat{R}&\rightarrow &  \widehat{\theta}_r\widehat{\theta}_r \\
\hspace{2cm}1 &\mapsto & \widehat{x}_r\otimes 1\otimes 1+1\otimes 1\otimes \widehat{x}_r.
\end{array}
\end{displaymath}
And finally for $m(s,r)=2$ we have 
\vspace{0.3cm}
\begin{displaymath}\label{d4}
\begin{array}{lll}
        \widehat{f}_{s,r} :\hspace{3cm}\widehat{\theta}_s\widehat{\theta}_r&\rightarrow &  \widehat{\theta}_r\widehat{\theta}_s \\
\widehat{R}\otimes_{\widehat{R}^s}\widehat{R}\otimes_{\widehat{R}^r}\widehat{R}\ni p_1\otimes p_2\otimes p_3 &\mapsto & p_1\widehat{\partial}_s(p_2)\otimes 1\otimes \widehat{x}_s p_3 +p_1\widehat{P}_s(p_2)\otimes 1\otimes p_3. 
\end{array}
\end{displaymath}

We define the following morphisms from  compositions of the previous ones :

\begin{itemize}\label{mor}
\item $\widehat{p}_r=(\mathrm{id}\otimes \widehat{j}_r)\circ (\widehat{\alpha}_r\otimes \mathrm{id}) :\widehat{\theta}_r\rightarrow \widehat{\theta}_r\widehat{\theta}_r ;$ 
\item $\widehat{\epsilon}_r=(\mathrm{id}\otimes \widehat{m}_r)\circ \widehat{\alpha}_r :\widehat{R} \rightarrow \widehat{\theta}_r ;$
\item $\widehat{x}_r= (m_r \circ \epsilon_r)/2 : \widehat{R} \rightarrow  \widehat{R}$
\end{itemize}

By calculating  with the explicit formulas we see that the morphism $\widehat{x}_r\in \mathrm{End}(\widehat{R})$ corresponds to the multiplication by $\widehat{x}_r\in \widehat{R}$.

\subsection{Graphical notation}

For sake of the clarity of the  exposition, we will introduce a graphical notation for the morphisms.  If we write $s_1\ s_2\cdots s_n$ in a diagram, it means $\widehat{\theta}_{s_1}\ \widehat{\theta}_{s_2} \cdots \widehat{\theta}_{s_n}$. 
 We identify $\widehat{R}$ with the blank space. Here we explain what symbol means what morphism: 
$$ \widehat{f}_{r,t}= \begin{array}{c}
\underbracket{r\ t}  \\
t\ r 
\end{array}
 $$

$$ \widehat{x}_s= {\diamond} $$

$$\widehat{\epsilon}_r= \begin{array}{c}
 \  \   \\
 r\hspace{-0.23cm}{^\smallfrown}\,     
\end{array}  \hspace{2cm}   \widehat{m}_r=\begin{array}{c}
r\hspace{-0.26cm}{_\smallsmile} \\
    
\end{array}  
 $$

$$\widehat{p}_r= \begin{array}{c}
 r  \\
 \overbrace{r\ r}   
\end{array}  \hspace{2cm}   \widehat{j}_r=\begin{array}{c}
\underbrace{r\ r} \\
 r
\end{array}  
 $$

$$\widehat{\alpha}_r=\begin{array}{c}
\   \\
    \overgroup{r\ r}
\end{array}\hspace{2cm} \widehat{m}_r\circ \widehat{j}_r=
\begin{array}{c}
\undergroup{r\ r}  \\

\end{array} 
$$
Pictures are to be understood up to multiplication by non-zero scalars.
If we write downwards a sequence of sequences of elements of $\mathcal{S}$, it means a chain of morphisms between the corresponding bimodules. For example if we write $\begin{array}{c}r\ s\ t \\t\ r\ s\\ r\ s\end{array}$ it means a chain of morphisms
 $$\widehat{\theta}_r\widehat{\theta}_s\widehat{\theta}_t \rightarrow \widehat{\theta}_t\widehat{\theta}_r\widehat{\theta}_s\rightarrow \widehat{\theta}_r\widehat{\theta}_s.$$
As an example, $\begin{array}{c}r\\\overbrace{\underbrace{r\ r}}\\ r
                                                                       \end{array}
$ means $\widehat{j}_r\circ \widehat{p}_r.$   We lose some information when we look at the picture of $\widehat{x}_s$, because we do not know by which $\widehat{x}_s$ are we multiplying, but this will not be a problem for our purposes.  So when we have a morphism we can draw its corresponding picture, but in general when we have a picture we cannot completely recover the morphism.

We will use bold letters when we mean a linear combination of this type of morphisms. For example  $\begin{array}{c}
 \ \mathbf{r}  \\ ^{\diamond}\mathbf{r}\\
 \overbrace{\mathbf{r}\ \mathbf{r}}\end{array} $ means  $\sum_i \lambda_i \widehat{p}_r\circ (\widehat{x}_{s_i}\otimes \mathrm{id}_{\widehat{\theta}_r})$, with $\lambda_i\in \widehat{R}$ and $s_i\in \mathcal{S}.$

\subsection{The tensor category $T_r$}\label{Tr}
For an introduction to tensor  categories and of strict tensor categories with a presentation by generators and relations, we refer the reader to \cite[chapters XI and XII]{Ka}.

\vspace{0.5cm}  \begin{definition}\label{d5} We  define the strict tensor category $T_r$ by generators and relations. Its objects are generated by $\theta_r$ and by the unit $R$ as tensor category. Its morphisms are generated by :
\begin{itemize}
 \item $j_r:\theta_r \rightarrow \theta_r\theta_r$
\item  $m_r:\theta_r\rightarrow R$
\item $\alpha_r: R \rightarrow \theta_r\theta_r $
\end{itemize}
 We define $p_r,\epsilon_r$ and $x_r$ by the same formulas as in section \ref{mor} without the hats. The relations defining $T_r$ are the following :
\begin{enumerate}
\item \label{6}$ \epsilon_r=(m_r\otimes \mathrm{id})\circ
 \alpha_r  \hspace{5.6cm}  \begin{array}{c} \overgroup{r\ r\hspace{-0.26cm}{_\smallsmile}}\\ r\ \ 
 
\end{array}=\begin{array}{c} \overgroup{r\hspace{-0.26cm}{_\smallsmile} \  r}\\ \ \ r
 
\end{array}   $
   \vspace{0.5cm} \item \label{7}$ p_r=(j_r\otimes \mathrm{id})\circ (\mathrm{id}\otimes \alpha_r)   \hspace{4cm} \begin{array}{c}\ \ \ \ \  r\\\overgroup{r\ \ r} \hspace{-0.45cm}\underbrace{\ \ r}\\ \hspace{-0.2cm}r\ \ \hspace{0.2cm}r
 
\end{array}=\begin{array}{c} r\ \ \ \ \  \\ \underbrace{r\ \ } \hspace{-0.45cm}\overgroup{r\ \ r} \\ \hspace{0.25cm}r\ \ \hspace{0.2cm}r

\end{array}$

\vspace{0.5cm} \item \label{1}$(\mathrm{id}\otimes (m_r\circ j_r))\circ ( \alpha_r \otimes
 \mathrm{id})=\mathrm{id} \hspace{3.1cm}  \begin{array}{c}
 \ \ \  \ \ r  \\
    \undergroup{r\;r} \hspace{-1.2cm}\overgroup{r\ \ }\\r\ \ \ \ \ \\
\end{array} = \begin{array}{c}r\\ \ \\ \ \end{array}$
\vspace{0.4cm}
\vspace{0.5cm} \item \label{y2}$((m_r\circ j_r)\otimes \mathrm{id})\circ (\mathrm{id}\otimes
 \alpha_r)=\mathrm{id} \hspace{2.8cm}  \begin{array}{c}
 r  \ \ \ \ \  \\
   \undergroup{r\;r}\hspace{-0.4cm}\overgroup{\ \ r}\\
\ \ \ \ \ r
\end{array} = \begin{array}{c}r\\ \ \\ \ \end{array}$
\vspace{0.5cm} \item \label{3} $ j_r\circ (\mathrm{id}\otimes j_r)=j_r\circ (j_r\otimes \mathrm{id})
\hspace{3.5cm}\begin{array}{c} r\underbrace{r\ r}\\ \hspace{-0.3cm}\underbrace{r\ \ r}\\ r 
\end{array}
  =\begin{array}{c}\underbrace{r\ r}r\\\hspace{0.3cm}\underbrace{r\ \ r}\\ r 
\end{array}$

    \vspace{0.5cm} \item \label{5}$  j_r\circ \alpha_r=0  \hspace{7.3cm} \begin{array}{c} \overgroup{\underbrace{rr}}
 
\end{array}=0 $
    
\vspace{0.5cm} \item \label{12}$     \mathrm{id}\otimes m_r=m_r\otimes \mathrm{id}+j_rx_r-x_rj_r   \hspace{2cm}  \begin{array}{c}  r\ r\hspace{-0.26cm}{_\smallsmile}\\r\ \  
\end{array}=\begin{array}{c} r\hspace{-0.26cm}{_\smallsmile}\ r\\ \ \ r 
\end{array}+\begin{array}{c}\underbrace{r\ r}\\ r\\ \ r^{\diamond}
\end{array}+\begin{array}{c}\underbrace{r\ r}\\ r\\  ^{\diamond}r\ \  
\end{array}$

    \vspace{0.5cm} \item \label{11}$  j_r\circ (\mathrm{id}\otimes x_r \otimes \mathrm{id})=$ 

\vspace{-0.6cm}
$=(m_r\otimes
 \mathrm{id})-(x_r\otimes \mathrm{id})\circ j_r  \hspace{3.2cm}   \begin{array}{c}r\ r\\\underbrace{r^{\diamond}r}\\r
 
\end{array}=\begin{array}{c}r\hspace{-0.26cm}{_\smallsmile}\ r\\ \ \ r
 
\end{array}+\begin{array}{c}\underbrace{r\ r}\\r\\^{\diamond}r\ \ 
 
\end{array}$

\end{enumerate}
\end{definition}

The following definition is needed to state the proposition \ref{tri}.

\vspace{0.5cm}  \begin{definition}\label{d6}
 Let $A$ be a commutative ring and $\mathcal{C}$ an $A-$linear category. If $R$ is a commutative $A-$algebra, we define the category $\mathcal{C}\otimes_{A}R$ in the following way : it has the same objects as $\mathcal{C}$ and its morphisms are defined by the formula 
$$\mathrm{Hom}_{\mathcal{C}\otimes_{A}R}(M,N)= \mathrm{Hom}_{\mathcal{C}}(M,N)\otimes_{A}R$$  
\end{definition}

Let $\mathcal{A}$ be the subring of $\mathrm{End}_{T_r}(R)$ generated  the set $\{x_s\}_{s\in \mathcal{S}}$.

\vspace{0.5cm}  \begin{proposition}\label{tri}
 Let $W$ be a Weyl group of type $A_1$ and let $\mathcal{A}$ be defined as above. The application $x_s\mapsto \widehat{x_s}$ extends to a morphism $\mathcal{A}\rightarrow \widehat{R}$ and there exists an equivalence of $\widehat{R}$-linear tensor categories between $\widehat{\mathbf{B}}(W)$ and $T_r\otimes_{\mathcal{A}}\widehat{R}$.
\end{proposition}

This proposition is a special case of theorem \ref{bacanz} that will be proved in the sequel.

\subsection{Some notation for the morphisms}\label{some}

Let $(s_1,\ldots,s_n)\in \mathcal{S}^n.$ In the article \cite{L}, some basis are constructed (called BFL basis) as a right $\widehat{R}-$module of $\mathrm{Hom}_{(\widehat{R},\widehat{R})}(\widehat{\theta}_{s_1}\cdots \widehat{\theta}_{s_n},\widehat{R}).$ We will fix one  such basis by taking (with the notations of the paper \cite{L}) $p(n,x,\overline{t})$ a $m-$tuple with $m$ minimal. To recall what is this basis, we start with a definition. 

\vspace{0.5cm}  \begin{definition}\label{h} 
If they have a meaning, we define in  $\mathrm{Hom}(\widehat{\theta}_{u_1}\cdots \widehat{\theta}_{u_l},\widehat{\theta}_{t_1}\cdots
 \widehat{\theta}_{t_k})$ the following morphisms,
\begin{enumerate}
    \item $^{i}\widehat{j}:= \mathrm{id}^i \otimes \widehat{j}_{u_{i+1}}\otimes \mathrm{id}^{l-i-2}$
    \item $^{i}\widehat{m}:=  \mathrm{id}^i  \otimes   \widehat{m}_{u_{i+1}}\otimes \mathrm{id}^{l-i-1}  
    $
    \item $^{i}\widehat{\alpha_s}:=  \mathrm{id}^i  \otimes  \widehat{\alpha}_r\otimes \mathrm{id}^{l-i}   $
    \item $^{i}\widehat{f}:= \mathrm{id}^i \otimes  \widehat{f}_{u_{i+1},u_{i+2}}\otimes
 \mathrm{id}^{l-i-2} $
    \item $^{i}\widehat{p}:= \mathrm{id}^i \otimes \widehat{p}_r\otimes \mathrm{id}^{l-i-1}$
 \item $^{i}\widehat{\epsilon_s}:= \mathrm{id}^i \otimes \widehat{\epsilon}_r\otimes \mathrm{id}^{l-i-1}$    
\item $^{i}\widehat{x_s}:= \mathrm{id}^i \otimes \widehat{x}_r\otimes \mathrm{id}^{l-i-1}$

Let $\mathrm{Mo}=\{j,m,\alpha_s,f,p,\epsilon_s, x_s\}$\label{d65}.
For $d\in$ Mo we define $\widehat{d}^{i}$ in almost the same way, the only difference is that we put $\mathrm{id}^i$ in the right. For example $\widehat{j}^{i}= \mathrm{id}^{l-i-2}\otimes \widehat{j}_{u_{i+1}}\otimes \mathrm{id}^i.$
\end{enumerate}
\end{definition}

Now we are able to define $\widehat{m}(t), \widehat{\mathrm{ch}}(t)$\label{d7} (chain) and $\widehat{\mathrm{cch}}(t) $ (complete chain) :\vspace{1cm}
\begin{itemize}
 \item $\widehat{m}(t)=\widehat{m}^{t} \hspace{8cm} r\hspace{-0.26cm}{_\smallsmile}$
\vspace{0.8cm}\item $\widehat{\mathrm{ch}}(t)=\widehat{j}^t\circ \widehat{f}^{t+1}\circ \widehat{f}^{t+2}\circ \cdots \circ \widehat{f}^{t'}   \hspace{0.5cm}\begin{array}{c} \underbracket{\ \ \ \  }\\\hspace{0.7cm}\underbracket{\ \ \ \ }\\\hspace{1.9cm}\ddots\\\hspace{3.1cm} \underbracket{\ \ \ \ }\\ \hspace{4.1cm}\underbrace{} \end{array}
$\vspace{1.5cm}
\item $\widehat{\mathrm{cch}}(t)=\widehat{m}^t\circ \widehat{\mathrm{ch}}(t)  \hspace{3cm}\begin{array}{c} \underbracket{\ \ \ \  }\\\hspace{0.7cm}\underbracket{\ \ \ \ }\\\hspace{1.9cm}\ddots\\\hspace{3.1cm} \underbracket{\ \ \ \ }\\ \hspace{4.1cm}\underbrace{}\\ \hspace{4.4cm} \hspace{-0.26cm}{_\smallsmile}\end{array}$
\end{itemize}

\vspace{0.5cm}  \begin{definition}\label{d8}
Let $\mathrm{No}= \{\widehat{m}, \widehat{\mathrm{ch}},\widehat{\mathrm{cch}}\}$. We say that $(\widehat{g_q}(t_q), \cdots , \widehat{g_1}(t_1))$, with $\widehat{g}_i\in      \mathrm{No}$, is a good $g$-expression of the morphism $\widehat{g_q}(t_q)\circ \cdots \circ \widehat{g_1}(t_1)$. If $t_{p+1}<t_p$ for all $1\leq p\leq q$ then we say that the good $g$-expression is in the good order.
\end{definition}

\vspace{0.5cm}  \begin{definition}\label{d9}
 Let $\nu=(\widehat{g_q}(t_q), \cdots, \widehat{g_1}(t_1))$ be a good $g$-expression. If $m\leq q$ and $\widehat{g_m}(t_m)\circ \cdots \circ \widehat{g_1}(t_1):\widehat{\theta}_{s_1}\cdots \widehat{\theta}_{s_n}\rightarrow \widehat{\theta}_{u_1}\cdots \widehat{\theta}_{u_l}$ we say that Ima$(\nu, m)=(u_1,\ldots, u_l).$
\end{definition}

\vspace{0.5cm}  \begin{definition}\label{d10}
 Let $\bar{t}=(t_1,\ldots, t_p)\in \mathcal{S}^p.$ We say that the $i^{th}$ element of $\bar{t}$ is of left type if there exists $j<i$, with $t_j=t_i$ and $m(t_r, t_i)=2$ for all $j< r<i.$
\end{definition}

\subsection{A basis of the morphism spaces in the right-angled case}\label{estria}
We start with a definition.
\vspace{0.5cm}  \begin{definition}
 We say that a good $g$-expression in good order $(\widehat{g_q}(t_q), \cdots, \widehat{g_1}(t_1))$ satisfies  property (P) if  for all $1\leq m\leq q$ such that the $t_m^{th}$ element of  Ima$(\nu, m)$ is of left type, we have  $\widehat{g_{m+1}}(t_{m+1})=\widehat{\mathrm{ch}}(t_m-1)$ or $\widehat{\mathrm{cch}}(t_m-1)$.
\end{definition}

The following is a corollary of theorem \cite[thm. 5.1]{L}  :

\vspace{0.5cm}  \begin{proposition}\label{FL}
Let $(W,\mathcal{S})$ be a right-angled Coxeter system. Let $(s_1,\ldots , s_n)\in \mathcal{S}^n$. The set 
\begin{multline*}
\widehat{FL}(s_1,\ldots,s_n)=\{  \ \mathrm{good\ }g\mathrm{-expressions\ in\ good\ order\  of\  morphisms\ in} \\ \mathrm{Hom}(\widehat{\theta}_{s_1}\cdots \widehat{\theta}_{s_n},\widehat{R} )  \mathrm{\ satisfying\  property}\ (P)\}. 
\end{multline*}
is an $\widehat{R}$-basis of $\mathrm{Hom}_{(\widehat{R},\widehat{R})}(\widehat{\theta}_{s_1}\cdots \widehat{\theta}_{s_n},\widehat{R} )$ called "Light leaves basis".
\end{proposition}

For every $M,N\in \widehat{\mathbf{B}}(W)$, we define the two morphisms :
\begin{itemize}\label{escroti}
    \item $\widehat{\textit{F}_s}(M,N) : \mathrm{Hom}(\widehat{\theta}_sM,N)\rightarrow \mathrm{Hom}(M,
 \widehat{\theta}_sN)$ that sends $f$ to $(\mathrm{id}_{\widehat{\theta}_s}\otimes f)\circ
 (\widehat{\alpha}_s\otimes \mathrm{id}_M)$
    \item $\widehat{\textit{G}_s}(M,N) : \mathrm{Hom}(M, \widehat{\theta}_sN)\rightarrow
 \mathrm{Hom}(\widehat{\theta}_sM,N)$ that sends $g$ to $((\widehat{m}_s\circ \widehat{j}_s)\otimes \mathrm{id}_N)\circ
 (\mathrm{id}_{\widehat{\theta}_s}\otimes g)$
\end{itemize}

They are inverse to each other (see \cite[lemma 3.3]{L} and its proof), so we can define a basis of  $\mathrm{Hom}(\widehat{\theta}_{s_1}\cdots \widehat{\theta}_{s_n},\widehat{\theta}_{t_1}\cdots \widehat{\theta}_{t_k})$ :
$$\widehat{FL}(s_1,\cdots,s_n ;t_1,\cdots,
 t_k):=\widehat{\textit{F}_{t_1}}\circ \cdots \circ \widehat{\textit{F}_{t_{k-1}}} \circ \widehat{\textit{F}_{t_k}}\circ \widehat{FL}(s_1, \cdots
 ,s_n).$$

\section{The tensor category \textbf{T}}\label{section3}

\subsection{}

\vspace{0.5cm}  \begin{definition}\label{d11} Let $(W,\mathcal{S})$ be a right-angled Coxeter system and $V$ a reflection faithful representation of $W$. We will define the tensor category $\mathbf{T}(W,V)$ by generators and relations. Its objects are generated as tensor category by $\theta_r$ for all $r\in\mathcal{S},$ and by the unit $R$. Its morphisms are generated by 
\begin{itemize}
 \item $j_r:\theta_r \rightarrow \theta_r\theta_r$
\item  $m_r:\theta_r\rightarrow R$
\item $\alpha_r : R \rightarrow \theta_r\theta_r $
\item $f_{sr} : \theta_s\theta_r\rightarrow \theta_r\theta_s$
\end{itemize}
We define $p_r,\epsilon_r$ and $x_r$ by the the same formulas as in section \ref{mor} without the hat.
For $s,r\in \mathcal{S}$, let us write (see section \ref{stab}) $$\widehat{P}_t(\widehat{x}_s)=\sum_{r\in \mathcal{S}}\lambda_{t,s}^r\widehat{x}_r \hspace{0.5cm} , \hspace{0.5cm} \widehat{I}_t(\widehat{x}_s)=\mu_{t,s} \widehat{x}_t \hspace{0.5cm} \mathrm{and} \hspace{0.5cm} \widehat{\partial}_t(\widehat{x}_s)=\frac{\widehat{I}_t(\widehat{x}_s)}{\widehat{x}_t}$$ 
We define in the same way $$P_t(x_s)=\sum_{r\in \mathcal{S}}\lambda_{t,s}^rx_r \hspace{0.5cm} , \hspace{0.5cm} I_t(x_s)=\mu_{t,s} x_t \hspace{0.5cm} \mathrm{and} \hspace{0.5cm} \partial_t(x_s)=\frac{I_t(x_s)}{x_t},$$ where $\lambda_{t,s}^r$ and $ \mu_{t,s}$ are elements of the field $k$.

 We define a monomial as a product of $x_s'$\label{d12}s, and a polynomial $\lambda \in \mathrm{End}(R)$ as a linear combination of monomials with coefficients in $k$. We will note $\mathcal{A}\subseteq \mathrm{End}(R)$ the subring of polynomials.

 The relations defining $T$ are the first 7 relations of definition  \ref{d5} and the following new relations (from \textbf{a} to \textbf{e} we assume that $m(s,r)=m(r,t)=m(s,t)=2$)  :
\begin{description}

\item[a)] \label{13}$ f_{sr}\circ f_{r,s}=\mathrm{id}   \hspace{7.6cm} \begin{array}{c}
 \underbracket{sr}  \\
 \underbracket{rs}\\
sr
\end{array} =\begin{array}{c}sr\\ \ \ \\ \ \ 
\end{array}$

\vspace{0.5cm}

    \item[b)] \label{b}$ (m_r\otimes \mathrm{id})\circ f_{sr}=\mathrm{id}\otimes m_r \hspace{5.4cm} \begin{array}{c}
 \underbracket{sr}  \\
 r\hspace{-0.26cm}{_\smallsmile}s\\
\ s
\end{array} =\begin{array}{c}s\,r\hspace{-0.26cm}{_\smallsmile}\\ s \ \,\\ \ \end{array}  $

\vspace{0.5cm}

    \item[c)] \label{16}$  (j_s\otimes \mathrm{id}) \circ (\mathrm{id}\otimes f_{r,s})=  
  f_{r,s} \circ (\mathrm{id}\otimes j_s) \circ (f_{sr}\otimes \mathrm{id})  \hspace{1cm} \begin{array}{c}
                                                                         s\, \underbracket{r\,\, \, s}\\ \underbrace{s\ s} r \\
\ \  s \hspace{0.4cm}r    \end{array}=\begin{array}{c} \underbracket{s\,\,\, r}\, s\\ \ r\underbrace{s\ s}\\\underbracket{r\ \ s}\\s\ \ r
\end{array}
$

\vspace{0.5cm}

\item[d)] $(\mathrm{id}\otimes f_{sr})\circ (\alpha_s\otimes \mathrm{id})=$

\vspace{-0.6cm}
$
\hspace{2cm}=(f_{r,s}\otimes \mathrm{id})\circ (\mathrm{id}\otimes \alpha_s)   \hspace{2.5cm} \begin{array}{c} \ \ \ \ \ r\\ \overgroup{s\ s} \hspace{-0.4cm}\underbracket{\ \ \ r}\\s\ r\ s\end{array}= \begin{array}{c} r\ \ \ \ \ \\  \underbracket{r\ s}\hspace{-0.3cm} \overgroup{\ \ s} \\ s\ r\ s\end{array}
 $

\vspace{0.5cm}

    \item[e)] \label{17}\vspace{0.4cm}$  (\mathrm{id}\otimes f_{sr}) \circ (f_{s,t}\otimes
 \mathrm{id}) \circ (\mathrm{id}\otimes f_{r,t}) = $\vspace*{-0.8cm}

\hspace*{1cm}$(f_{r,t}\otimes \mathrm{id})  \circ (\mathrm{id}\otimes f_{s,t})
  \circ (f_{sr}\otimes \mathrm{id}) $\vspace{0.5cm} $\hspace*{2cm}  \begin{array}{c}
                                                              s\underbracket{r\ t}\\ \underbracket{s\ t}r\\ t \underbracket{s\ r}\\t\ r\, s
                                                             \end{array}= \begin{array}{c}
                                                              \underbracket{s\ r}t\\ r\underbracket{s\ t}\\ \underbracket{r\ t}s\\t\ r\, s
                                                             \end{array}$
    
\vspace{0.5cm}

\item[f)]\label{kra}  $f_{t,r}\circ (\mathrm{id}_{\theta_r}\otimes x_r\otimes \mathrm{id}_{\theta_r})=$

\vspace{-0.6cm}$ (P_t(x_s)\otimes \mathrm{id}_{\theta_r\theta_t}) \circ f_{t,r}+(\mathrm{id}_{\theta_r\theta_t}\otimes I_t(x_s))\circ f_{t,r}  \hspace{1.3cm} \begin{array}{c}t\ r\\ \underbracket{t^{\diamond}r}\\ r\ t
                                                                                                \end{array} = \begin{array}{c} \underbracket{\mathbf{t}\ \mathbf{r}}\\\mathbf{r}\ \mathbf{t} \\ ^{\diamond}\mathbf{r}\ \mathbf{t}\ 
                                                                                                \end{array}  + \begin{array}{c} \underbracket{\mathbf{t\ r}}\\\mathbf{r} \ \mathbf{t} \\\ \mathbf{r}\ \mathbf{t}^{\diamond}
                                                                                                \end{array}
 $

\vspace{0.7cm}

\item[g)]  $j_r\circ (\mathrm{id}_{\theta_r}\otimes x_r\otimes \mathrm{id}_{\theta_r})=\partial_r(x_s)m_r\otimes \mathrm{id}+$

\vspace{-0.7cm}\hspace{0.5cm}$+(P_r(x_s)\otimes \mathrm{id}- x_s\partial_r(x_s)\otimes \mathrm{id})\circ j_r   \hspace{2.2cm} \begin{array}{c}r\ r\\ \underbrace{r^{\diamond}r}\\ r
                                                                                                \end{array} = \begin{array}{c}r\hspace{-0.26cm}{_\smallsmile}\ r \\ \ r
                                                                                                \end{array}    +\begin{array}{c} \underbrace{\mathbf{r}\ \mathbf{r}}\\ \mathbf{r}\\ ^{\diamond}\mathbf{r}\ 
                                                                                                \end{array} 
 $

\vspace{0.3cm}

\end{description}

\end{definition}

\vspace{0.5cm}  \begin{remark}
 We can say heuristically that \textbf{a}, \textbf{b}, \textbf{c} and \textbf{d} say that the $f_{sr}$ "commute" with everything, \textbf{e} is the hexagon relation typical in the symmetric monoidal categories. The last two relations tells us how to "take down" the $x_s'$s.
\end{remark}

 The ring $\mathrm{End}(R)$ is commutative because of general results in tensor categories (see \cite[prop. X1.2.4]{Ka}). As before, for $s\in\mathcal{S}$, we make $x_s\in \mathrm{Arr}(\mathbf{T}(W,V))$ to act in $\widehat{R}$ in the same way as $\widehat{x_s}\in \widehat{R}$ does. We deduce an action of $\mathcal{A}$ in $\widehat{R}$.

The following is the central result of this paper :

\vspace{0.5cm}  \begin{theorem}\label{bacanz}
Let $(W,\mathcal{S})$ be a right-angled Coxeter group. The functor that sends $\theta_s$ to $\widehat{\theta_s}$,  $j_s\otimes 1$ to $\widehat{j_s}$, etc. is an  equivalence of $\widehat{R}$-linear tensor categories between $\mathbf{T}(W,V)\otimes_{\mathcal{A}}\widehat{R}$ and  $\widehat{\mathbf{B}}(W,V)$.
\end{theorem}

\textbf{Proof.}
   We  repeat all definitions in sections  \ref{some} and \ref{estria} taking out the hat, so we define $^{i}d$ and $ d^{i},$ for $d\in$Mo, $ m(t), \mathrm{ch}(t), \mathrm{cch}(t)$, we define a good $g$-expression and to be of left type, $FL$, $F_s(M,N)$ and $G_s(M,N)$ 

We will put $\mathbf{T}=\mathbf{T}(W,V)\otimes_{\mathcal{A}}R$ and $\widehat{\mathbf{B}}=\widehat{\mathbf{B}}(W,V)$.
We define a  functor $\mathfrak{Fu}$\label{d13} from $\mathbf{T}$   to $\widehat{\mathbf{B}}$.  We define $\mathfrak{Fu}(R)=\widehat{R}$, and for all $s\in\mathcal{S}$, $\mathfrak{Fu}(\theta_s)=\widehat{\theta}_s$. If $M$ and $M'$ are objects of $\mathbf{T}$, we define $\mathfrak{Fu}(M\otimes M')=\mathfrak{Fu}(M)\otimes \mathfrak{Fu}(M').$

For all $s\in\mathcal{S}$ we define $\mathfrak{Fu}(j_s\otimes \lambda)=\widehat{j}_s\lambda$, $\mathfrak{Fu}(m_s\otimes \lambda)=\widehat{m}_s\lambda$, etc. As we know explicitly all the morphisms in $\widehat{\mathbf{B}}$, we can  easily verify that all the  relations are satisfied in $\widehat{\mathbf{B}}$, so by  \cite[proposition XII.1.4]{Ka} we have that $\mathfrak{Fu}$ defines a tensor functor. By  \cite[thm. 5.1]{L} we know that the set of morphisms $\{\widehat{j}_r,\widehat{m}_r,\widehat{\alpha}_r,\widehat{f}_{s,r}\}$ generates (as tensor category) all the morphisms in $\widehat{\mathbf{B}}$. So we only need to prove that for all $M,N \in \mathbf{T}$, 
the map $\mathfrak{Fu} : \mathrm{Hom}_{\mathbf{T}}(M,N) \rightarrow \mathrm{Hom}_{\widehat{\mathbf{B}}}(\mathfrak{Fu}(M),\mathfrak{Fu}(N))$ is injective.

 We  start by proving the following 
\begin{lemma}
The applications $\textit{F}_s(M,N)$ and $\textit{G}_s(M,N)$ are  inverse to each other  
\end{lemma}
\textbf{Proof. } \begin{displaymath}
\begin{array}{lll}
 \textit{F}_s(M,N)\circ
 \textit{G}_s(M,N)(g)&=&\textit{F}_s(M,N)(((m_s\circ j_s)\otimes \mathrm{id}_N)\circ (\mathrm{id}_{\theta_s}\otimes g)) \\
 &=& \left\lbrace \mathrm{id}_{\theta_s}\otimes \left[((m_s\circ j_s)\otimes
 \mathrm{id}_N)\circ (\mathrm{id}_{\theta_s} \otimes g)\right]\right\rbrace \circ
 (\alpha_s\otimes \mathrm{id}_M)         \\
&=& \left[ (\mathrm{id}_{\theta_s}\otimes(m_s\circ j_s)\otimes \mathrm{id}_N)\circ
 (\mathrm{id}^2_{\theta_s}\otimes g)  \right]\circ(\alpha_s\otimes \mathrm{id}_M)\\
&=& (\mathrm{id}_{\theta_s}\otimes(m_s\circ j_s)\otimes \mathrm{id}_N)\circ
 (\alpha_s\otimes g)\\
&=& (\mathrm{id}_{\theta_s}\otimes(m_s\circ j_s)\otimes \mathrm{id}_N)\circ
 (\alpha_s\otimes \mathrm{id}_{\theta_s}\otimes \mathrm{id}_N)\circ g \\
&=& (\left\{ \left[  \mathrm{id}_{\theta_s}\otimes (m_s\circ j_s) \right]\circ
 (\alpha_s\otimes \mathrm{id}_{\theta_s})     \right\}\otimes \mathrm{id}_N)\circ g\\
&=& (\mathrm{id}_{\theta_s}\otimes \mathrm{id}_N)\circ g\\
&=&g
\end{array}
\end{displaymath}

If (0) is the relation   $(f_1\otimes g_1)\circ (f_2\otimes
 g_2)=(f_1\circ f_2)\otimes (g_1\circ g_2)$ (relation satisfied in all tensor categories), then all but the next to last equality  are derived from (0).  The next to last equality is derived from relation \ref{1} (see definition \ref{d5}). In a similar way we prove that  $
 \textit{G}_s(M,N)\circ\textit{F}_s(M,N)(f)=f$ using relation \ref{y2}. $\hfill  \Box$
\vspace{0.5cm}  \begin{definition}
 Consider two morphisms $g :M_1 \rightarrow N_1$ and $f :M_2\rightarrow N_2.$ The relation $(f\otimes \mathrm{id}_{N_1})\circ (\mathrm{id}_{M_2}\otimes g)= (\mathrm{id}_{N_2}\otimes g)\circ (f\otimes \mathrm{id}_{m_1})$, satisfied in all tensor categories, will be called commutation relation.
\end{definition}

\subsection{}
\vspace{0.5cm}  \begin{lemma}\label{bla}
Let us suppose that for every sequence $(s_1,\cdots,s_n)$, $f\in \mathrm{Hom}(\theta_{s_1}\cdots \theta_{s_n},R)$ implies that  $f=\sum_i a_i\lambda_i$, where $a_i\in FL(s_1,\cdots,s_n)$ and the $\lambda_i$ are polynomials. Then  for every couple of sequences  $(s_1,\cdots,s_n)$ and $(t_1,\cdots,t_k)$, $
 f'\in \mathrm{Hom}(\theta_{s_1}\cdots \theta_{s_n},\theta_{t_1}\cdots \theta_{t_k})$
 implies that  $f'=\sum_i a'_i\lambda'_i$, where $a'_i\in
 FL(s_1,\cdots,s_n;t_1\cdots t_k)$ and the $\lambda'_i$ are polynomials.
\end{lemma}
\textbf{Proof. } By hypothesis, $\textit{G}_{t_k}\circ \cdots \circ
 \textit{G}_{t_1} (f') =\sum_i a_i\lambda_i$, where $a_i\in FL(t_k,\cdots,
 t_1,s_1,\cdots,s_n)$ and the $\lambda_i$ are polynomials.
Then

\begin{displaymath}
\begin{array}{lll}
f'&=&\textit{F}_{t_1}\circ \cdots \circ \textit{F}_{t_k}\circ
 \textit{G}_{t_k}\circ \cdots \circ \textit{G}_{t_1} (f') \\
&=& \sum_i\textit{F}_{t_1}\circ \cdots \circ \textit{F}_{t_k}
  (a_i\lambda_i)  \\
&=&\left\{\sum_i\textit{F}_{t_1}\circ \cdots \circ \textit{F}_{t_k}
  (a_i)\right\}\lambda_i
\end{array}
\end{displaymath}

and by definition $\sum_i\textit{F}_{t_1}\circ \cdots \circ
 \textit{F}_{t_k}  (a_i)\in FL(s_1,\cdots,s_n;t_1\cdots t_k)$. $\hfill \Box$

\subsection{}\label{rest}From this lemma we can conclude that  to complete the proof of theorem \ref{bacanz}, it suffices to prove that for every sequence $(s_1,\cdots,s_n)$, $f\in \mathrm{Hom}(\theta_{s_1}\cdots \theta_{s_n},R)$ implies that  $f=\sum_i a_i\lambda_i$, where $a_i\in FL(s_1,\cdots,s_n)$ and the $\lambda_i$ are polynomials.

\vspace{0.5cm}  \begin{definition}\label{d14}
For a sequence $s_1,\ldots , s_n \in \mathcal{S}$, 
we will say that  $\theta_{s_1}\cdots \theta_{s_n}$ is a basic bimodule.
\end{definition}

Every morphism between basic bimodules can be written as a sum of terms of the following type :
 $^{i_\omega}\!d_{\omega}\circ \cdots \circ ^{i_1}\!\!d_1 $, with $d_i\in \mathrm{Lo}=\{j_s,m_s,f_{sr},\alpha_s\ \vert\ s, r\in \mathcal{S}\}$.

\vspace{0.5cm}  \begin{definition}\label{d15}In the sequel an \textit{expression} of a morphism $g$ between basic bimodules means a sequence $(^{i_\omega}\!d_{\omega}, \ldots , ^{i_1}\!\!d_1) $, with $d_i\in \mathrm{Lo}$, such that $^{i_\omega}\!d_{\omega}\circ \cdots \circ ^{i_1}\!\!d_1 =g.$ If $\nu$ is an expression of $g$, we define $\overline{\nu}=g,$\label{d18} and we say that $\nu$ is representing $g.$ An \textit{R-expression} is an expression representing a morphism from a basic bimodule to $R$. Sometimes we will consider a good $g$-expression simply like an expression in the obvious way.  Finally, if we consider a formal $\mathcal{A}$-linear combination of expressions $\nu=\sum_{i\in I}\lambda_i\nu_i$, with $\lambda_i\in \mathcal{A}$, we define $\overline{\nu}=\sum_{i\in I}\lambda_i\overline{\nu_i}$ 
\end{definition}

\vspace{0.5cm}  \begin{definition}\label{d16}If
  $(^{i_\omega}\!d_{\omega}, \ldots, ^{i_1}\!\!d_1) $ is an  expression of $g$, we say that the   $k^{th}$ term of this  expression is the morphism
 $^{i_k}\!d_{k}$. We say that 
 $^{i_\omega}\!d_{\omega}$ is the last term of this expression,
 $^{i_1}\!d_{1} $ is the first term, and $\omega$ is the length of the expression. 
\end{definition}

 When $i$ is not important, we will write $d$ instead of $^id$, and when the $s$ is not important we will write for example $x$ for $x_s$ and $f$ for $f_{sr}$. 

\vspace{0.5cm}  \begin{definition}\label{d19}
We define $G$ as the set of all $R$-expressions $(^{i_\omega}\!d_{\omega}, \cdots, ^{i_1}\!\!d_1) $ with $d_i\neq \alpha_s$ for all $i$ and for all $s\in\mathcal{S}$. 
\end{definition}

\vspace{0.5cm}  \begin{definition}
We will say that $\nu= (^{i_\omega}\!h_{\omega},\ldots, ^{i_1}\!\!h_1)$ with $h_i\in$ Mo for all $1\leq i\leq \omega$ is a generalized expression of $^{i_\omega}\!h_{\omega}\circ \cdots \circ ^{i_1}\!\!h_1$. 

 Let $$^{i_\omega}\!h_{\omega}\circ \cdots \circ ^{i_1}\!\!h_1=\sum_{z} \ ^{p_{u,z}}\ \!\!h'_{u,z}\circ \cdots \circ ^{p_{1,z}}\!h'_{1,z}$$ be one of the relations 1-8 or a-g defining the tensorial category \textbf{T}, where all $h_i$ and $h'_{i,z}$ are elements of Mo. For all natural numbers $b\geq 0$   we say that $X=(^{i_\omega+b}h_{\omega}, \ldots, ^{i_1+b}\!h_1)$ is a left generalized expression and $\sum_{z} (^{p_{u,z}+b} h'_{u,z}, \ldots, ^{p_{1,z}+b}\!h'_{1,z})$  is a right generalized expression correspondig to $X$.
\end{definition}

\vspace{0.5cm}  \begin{definition}\label{d20}
Let us consider  $\nu$ a generalized expression. We will say that we \emph{apply} to $\nu$ one of the relations 1-8 or a-g, defining the tensorial category \textbf{T}, and we obtain a finite sum $\sum_i\nu'_i$ of generalized expressions if $\sum_i\nu'_i$ is the result of changing in $\nu$ a substring $X$ that is a left generalized expression by a right generalized expression corresponding to $X$.

\end{definition}

 \section{A brief summary of what we will do}\label{section4}

If we have an expression $\nu\in G$, we will start by defining $\nu\{\#m-\mathrm{bad}\}$ and $\nu\{\#j-\mathrm{bad}\} \in \mathbb{N},$ that are measures of how far is $\nu$ from being in the light leaves basis.  We will see that an expression $\nu$ is in $FL(s_1,\ldots, s_n)$ if and only if $\nu$ is a good $g$-expression in the good order, and $\nu\{\#m-\mathrm{bad}\}=\nu\{\#j-\mathrm{bad}\}=0$. 

In each section we start with an expression satisfying some properties. Then we apply some relations to this expression and we arrive to another expression (eventually an $\mathcal{A}$-linear combination of expressions) which represents the same morphism, but satisfying new properties. In the following list, we show  what properties satisfies the expression with which we start in each section, and what properties satisfies the expression to which we arrive : 
\begin{itemize}
\vspace{0.2cm} \item \textit{Section 4 }: $\nu\in G \leadsto \nu\in G,\  \nu\{\#m-\mathrm{bad}\}=0$ 
\vspace{0.2cm}\item  \textit{Section 5 }: $\nu\in G \leadsto \nu\in G,\  \nu$ a good $g$-expression 
\vspace{0.2cm}\item    \textit{Section 6} : $\nu \leadsto \nu\in G$  
\vspace{0.2cm}\item \textit{Section 7} : $\nu$ a good $g$-expression $\leadsto \nu$ a good $g$-expression in the good order  
\vspace{0.2cm}\item \textit{Section 8.1 }: $\nu\in G\leadsto \nu$ a good $g$-expression in the good order with $\nu\{\#m-\mathrm{bad}\}=0$   
\vspace{0.2cm}\item \textit{Section 8.2} : $\nu$ a good $g$-expression in the good order with $\nu\{\#m-\mathrm{bad}\}=0 \leadsto$ $\nu$ a good $g$-expression in the good order with $\nu\{\#m-\mathrm{bad}\}=\nu\{\#j-\mathrm{bad}\}=0$
\end{itemize}

\section{Some numbers associated to an expression}
\vspace{0.5cm}  \begin{definition}\label{d21}\label{d24}
 Let $\nu =(^{i_k}d_{k}, \cdots, ^{i_1}\!d_1)$ be an expression for $g$, and let the sequence of integers $(i'_1, i'_2,\ldots, i'_k)$ be such that  $\overline{\nu}=d^{\, i'_k}_{k}\circ \cdots \circ d^{\, i'_1}_{1}$. We define $\nu_r=(^{i_r}d_{r}, \cdots, ^{i_1}\!d_1)$ the truncation  of $\nu$ at $r$.\newline

If $\overline{\nu_r} : \theta_{s_1}\cdots \theta_{s_n}\rightarrow \theta_{t_1}\cdots \theta_{t_p},$ we define $\widehat{\nu_r}=(t_1,\ldots, t_p)\in \mathcal{S}^p.$\newline
For $b \in\{ j,m,\alpha,f \}$ we define $\nu[b]=\{p; 1\leq p \leq k,  d_p=b\}$.

\vspace{0.2cm} 
 We say that $r$ is $m$-bad (resp. $j$-bad) for $\nu$ if $d_r=m$ (resp. $d_r=j$) and the $(i_r+1)^{th}$ element of $\widehat{\nu_r}$ is of left type.

 We define the following sets associated to $\nu$ :
\begin{itemize}
 \vspace{0.2cm}\item $A_m(\nu)=\{r\ ;\ r\ $ is $m-$bad for $\nu  \}$ 
\vspace{0.2cm}\item $A_j(\nu)=\{r\ ;\ r\ $ is $j-$bad for $\nu  \}$
\end{itemize}

\vspace{0.2cm}We define the following elements associated with $\nu :$
\begin{itemize}
 \vspace{0.2cm}\item $\nce=$card$(A_m(\nu))\in \mathbb{N}$
\vspace{0.2cm}\item $\nun=$card$(A_j(\nu))\in \mathbb{N}$
\end{itemize}
\vspace{0.2cm}Let $\nu[j]=\{a_1,\ldots, a_p\}$ with $a_1<a_2<\ldots <a_p.$
\begin{itemize}
\vspace{0.4cm}\item $\ndo =(a_1+i'_{a_1},\ldots, a_k+i'_{a_k})\in \mathbb{N}^k$
\vspace{0.4cm}\item $\ntr=$card$(\nu[f])\in \mathbb{N}$
\vspace{0.4cm}\item $\ncu=\sum_{p\in \nu[f]}(i_p)  \in \mathbb{N}$
\vspace{0.4cm}\item $\nci=\sum_{p\in \nu[m] \cup \nu[j]}(p) \in \mathbb{N}$
\vspace{0.4cm}\item $\nse= \sum_{p\in \nu[m]}(k-p) \in \mathbb{N}$
\vspace{0.4cm}\item $
\nsi= \begin{cases} \mathrm{min}A_m(\nu)\,\,\,\, \text{ if } \,\, A_m(\nu)\neq \emptyset \\ \hspace{0.8cm}  0
\hspace{1.3cm}
 \text{ if } \,\, A_m(\nu) = \emptyset
\end{cases}
$
\vspace{0.4cm}\item $\noc=\mathrm{card}(\{p>\nsi \ \vert\  i_p=1\ \, \mathrm{ or } \ \,i_p=2\})\in \mathbb{N}$
\vspace{0.4cm}\item  $ \nnu= (\noc,i_{\nsi})\in \mathbb{N}^2$
\vspace{0.4cm}\item $
\ndi= \begin{cases} \mathrm{max}A_j(\nu)\,\,\,\, \text{ if } \,\, A_j(\nu)\neq \emptyset \\ \hspace{0.8cm}  0
\hspace{1.3cm}
 \text{ if } \,\, A_j(\nu) = \emptyset
\end{cases}
$\end{itemize}
\vspace{0.4cm}If $p\in \nu[1]$ (resp. $p\in \nu[2]$) $ \hat{\nu_p}=(t_1,\ldots,t_k)$ then we define $\pi_1(p)$ (resp. $\pi_2(p)$) $=\mathrm{card}\{m\leq i_p \ \vert \ t_m=t_{i_{p+1}}\}$. Finally we define
 $$\non= \sum_{p\in \nu[1]\cup \nu[2]}(\pi_1(p) +\pi_2(p))\in \mathbb{N}$$

\end{definition}

\section{Restriction to the case $\nce=0$}
\vspace{0.5cm}  \begin{proposition}\label{primera}
 For every $\nu \in G$ there exists a set $\Delta$, polynomials $\lambda_{\delta}$ and elements $\nu_{\delta}\in G$ such that $\overline{\nu}=\overline{\sum_{\delta\in \Delta}\lambda_{\delta}\nu_{\delta}}$ and such that $\nu_{\delta}\ce=0$ for all $\delta \in \Delta.$
\end{proposition}

\textbf{Proof. }

We will start by constructing  $\mathcal{F}_1(\nu), \mathcal{F}_2(\nu)$\label{d26}, linear combinations of expressions such that 
$\overline{\mathcal{F}_i(\nu)}=\overline{\nu}$ for $1\leq i \leq 2.$

\subsubsection{$\mathcal{F}_1(\nu)$ : ``Taking out the $x_r$'s''}\label{tak} \vspace{0.5cm}  \begin{definition}\label{d25}If $\nu=(^{i_p}d_p, \cdots,  ^{i_{r+1}}d_{r+1}, ^{i}x_s, ^{i_r}d_r, \cdots, ^{i_1}d_1)$ is a generalized expression, with $d_l\in \mathrm{Lo}-\{\alpha_s\}_{s\in \mathcal{S}}$ for all $1\leq l \leq p$, we say that $\nu \in G^p_r.$\end{definition}

Let $\nu \in G^p_r$. We take $x_s$ as far to the left as possible (or down in the picture) with commutation relations. If it does not arrive to the last term, we obtain a string of the form $\begin{array}{c}s\ s\\ \underbrace{s^{\diamond}s}\\ s \end{array}$ or $\begin{array}{c}t\ r\\ \underbracket{t^{\diamond}r}\\ r\ t\end{array}$. If we are in the first case we apply relation \textbf{g} and if we are in the second case, we apply relation \textbf{f} and we obtain $$\mathcal{F}_1^1(\nu)=\sum_{i\in I_1} w_i^1+ \sum_{b\in B_1}g_b^1,$$
with $\overline{\nu}= \overline{\mathcal{F}_1^1(\nu)}, $ and
where  $I_1$ and $B_1$ are finite sets, $g^1_b\in G$ (corresponding to $\begin{array}{c}
                                                                         r\hspace{-0.26cm}{_\smallsmile} \ r\\ \ \ r
                                                                        \end{array}
$ in relation \textbf{g}) and
\begin{equation}\label{trom}
 w_i^1 \in G^p_{r_i^1},\  \mathrm{with }\ \  r^1_i\geq r+1
\end{equation}

 We define $\mathcal{F}_1^n$ inductively : let  $\mathcal{F}_1^n(\nu)=\sum_{i\in I_n} w_i^n+ \sum_{b\in B_n}g_b^n.$ We define
$\mathcal{F}_1^{n+1}(\nu)=\sum_{i\in I_n} \mathcal{F}_1^1(w_i^n)+ \sum_{b\in B_n}g_b^n$. By induction and inequality (\ref{trom}) we deduce that
$w_i^N \in G^p_{r_i^N}$, with $r^N_i\geq r+N$. We conclude that for $N\gg 0$, $\mathcal{F}_1^N(\nu)=\mathcal{F}_1^{N+1}(\nu)$. We define $\mathcal{F}_1(\nu)=\mathcal{F}_1^{N}(\nu).$

As $\mathcal{F}_1(\nu)=\sum_{i\in I_N} w_i^N+ \sum_{b\in B_N}g_b^N$ with $w_i^N\in G^p_{p+1}$, and as $w\in G^p_{p+1}\Rightarrow w=\lambda g$ with $\lambda$ a polynomial and $g\in G$, we have :
\vspace{0.5cm}  \begin{lemma}
There exists polynomials $\lambda_i$ and elements $g_i\in G$ such that
 $\mathcal{F}_1(\nu)=\sum_{i\in I}\lambda_ig_i,$
\end{lemma}

\vspace{0.5cm}  \begin{remark}
 Essentially when we take the $x_r'$s out of $\nu$ the only thing  we change in the picture of $\nu$ is eventually changing some $j'$s by $m'$s.
\end{remark}

\subsection{$\mathcal{F}_2(\nu)$\label{d27}}

Let $\nu \in G$. If $\nsi=0$ we define $\mathcal{F}_2(\nu)=\nu.$ Now let us suppose $\nsi\neq 0.$ By definition, in the $\nsi^{th}$ position of $\nu$ we have a picture like this one : $\cdots r \cdots r\hspace{-0.26cm}{_\smallsmile}$, where $r$ commutes with all the elements between the two $r'$s. By relation \textbf{a.} applied several times we obtain 
$$ \cdots r \cdots r\hspace{-0.26cm}{_\smallsmile}=\begin{array}{c} \cdots \underbracket{\ r \ \ \  }\cdots \hspace{0.9cm}r \cdots \\ \hspace{0.4cm}\cdots \underbracket{\ r \ \ \  }\cdots \hspace{0.5cm}r \cdots\\\hspace{0.5cm}\ddots \hspace{0.5cm}\\\hspace{1.3cm} \underbracket{\ r \ \ \  }\ r\\ \hspace{2.5cm}r\ \  r\hspace{-0.26cm}{_\smallsmile}\ \cdots  
\\ \hspace{0.9cm}\cdots \underbracket{\ \ \  r\,  }\cdots
\\ \hspace{0cm} \udots
\\ \hspace{-1.2cm}\cdots \underbracket{\ \ \  r\,  }\cdots
\\ \hspace{-2.1cm}\cdots \underbracket{\ \ \  r\,}\cdots
   \end{array}$$

Or in formulas : $^ym=^x\hspace{-0.17cm}f\circ ^x\hspace{-0.17cm}f \circ ^{x+1}\hspace{-0.1cm}f\circ \cdots \circ ^{y-2}\hspace{-0.1cm}f\circ ^y\hspace{-0.17cm}m\circ ^{y-2}\hspace{-0.17cm}f\circ \cdots \circ ^{x+1}\hspace{-0.17cm}f\circ ^x\hspace{-0.17cm}f.$ In the right hand side of the picture we apply relation \ref{12}. We obtain $\overline{\nu}=\overline{\nu_1+\nu_2+\nu_3}$, where $\nu_1\in G$ and $\nu_2,\nu_3\in G^p_r$ for some $r$. We define $\mathcal{F}^1_2(\nu)=\nu_1+\mathcal{F}_1(\nu_2)+\mathcal{F}_1(\nu_3).$

If we write $\mathcal{F}^1_2(\nu)= \sum_i \lambda_i^1 g_i^1$, as 
\begin{multline*}\mathrm{card}(\{1\leq p\leq \nsi\ \vert\  i_p=1 \ \mathrm{or} \  i_p=2\})+\\+\noc= \mathrm{card}(\{p\ \vert\ i_p=1 \ \mathrm{or} \  i_p=2\}) \end{multline*}
is a constant for all the expressions of $\overline{\nu}$,  by construction we have 
\begin{equation}\label{nu9}
 \mathcal{F}^1_2(g_i)\nue<\nnu.
\end{equation}

If $\mathcal{F}^N_2(\nu)= \sum_i \lambda_i^N g_i^N$, we define $\mathcal{F}^{N+1}_2(\nu)= \sum_i \lambda_i^N\mathcal{F}_2^1(g_i^N)$, and by equation \ref{nu9} and the fact that $\nnu\geq (0,0) $ we conclude that there exists $N\gg 0$ such that $\mathcal{F}^N_2(\nu)=\mathcal{F}^{N+1}_2(\nu).$ We then define $\mathcal{F}^N_2(\nu)=\mathcal{F}_2(\nu).$ This means in particular that if 
$\mathcal{F}_2(\nu)=\sum_{\delta\in \Delta}\lambda_{\delta}\nu_{\delta},$ then $\nu_{\delta}\ce=0$, which proves proposition \ref{primera}
$\hfill \Box$

\section{}

\vspace{0.5cm}  \begin{proposition}\label{duru}For every $\nu \in G$ there exists a set $\Delta$, polynomials $\lambda_{\delta}$ and elements $\nu_{\delta}\in G$ such that $\overline{\nu}=\overline{\sum_{\delta\in \Delta}\lambda_{\delta}\nu_{\delta}}$ and such that 
\begin{itemize}
 \item $\nu_{\delta}\ce=0$ for all $\delta \in \Delta.$
\item If the $k^{th}$ term of $\nu_{\delta}$ is $^if$ then the $(k+1)^{th}$ term is $^{i+1}f$ or $^{i+1}j.$
\end{itemize}
\end{proposition}
\vspace{0.5cm}  \begin{remark}\label{ya}
 The second condition means that $\nu_{\delta}$ is a good $g$-expression.
\end{remark}

\subsection{}

\textbf{Proof. }

We start with some definitions.

\vspace{0.5cm}  \begin{definition}
\begin{itemize}
 \item Relation $\begin{array}{c} s\ \underbracket{s\ r}\\ \underbracket{s\ r}\ s\\r\underbrace{s\ s}\\\hspace{-0.4cm}r\ \ \ s \end{array}= \begin{array}{c} \underbrace{s\ s}r\\\underbracket{s\ \ r}\\r\ \ s \end{array}$ can be deduced from \textbf{a.} and \textbf{c.} and will be called \textbf{c'.}
\item Relation $ \begin{array}{c} \underbracket{r\ s}\\ s\  r\hspace{-0.26cm}{_\smallsmile}\\s\ \ \end{array}=\begin{array}{c} r\hspace{-0.26cm}{_\smallsmile}\ s\\ \ \ s\end{array}
 $ can be deduced from \textbf{a.} and \textbf{b.} and will be called \textbf{b'.}
\item A commutation relation of type $\begin{array}{c} f\\ j \end{array}= \begin{array}{c} j\\f \end{array}$ or 
$\begin{array}{c} f\\m \end{array} = \begin{array}{c} m\\f \end{array}$ will be called a relation \textbf{x.}
\item A commutation relation of type $\begin{array}{c} ^if\\^jf \end{array} = \begin{array}{c} ^jf\\^if \end{array}$ will be called relation \textbf{y.}
\end{itemize}
\end{definition}

\subsection{} 

\begin{definition}
 For $\nu\in G$ we define $Y_{\nu}$ to be the set of  $\nu'\in G$  such that there exists a natural number  $n$ and a sequence $(\nu_n,\ldots,\nu_2,\nu_1)$ satisfying that $\nu_1=\nu,$ $\nu_n=\nu'$ and $\nu_{i+1}$ is obtained by applying a relation \textbf{y.} to $\nu_i$ for all $1\leq i\leq n-1.$ 

We say that $\nu$ satisfies property (Q) if for all $\nu'\in Y_{\nu}$  the relations \textbf{a.}, \textbf{b.}, \textbf{b'.}, \textbf{c.}, \textbf{c'.}, \textbf{d.} or \textbf{x.} cannot be applied to $\nu'$. In other words, if $\nu'\in Y_{\nu}$, there is no generalized left expression in $\nu'$ corresponding to relations \textbf{a.}, \textbf{b.}, \textbf{b'.}, \textbf{c.}, \textbf{c'.}, \textbf{d.} or \textbf{x.}. 
\end{definition}

For $\nu\in G$ we will define $\mathcal{F}_3(\nu)$\label{d28}.
\vspace{0.5cm}  \begin{lemma}\label{5.3}
Let $\nu\in G$. There exists an expression $\mathcal{F}_3(\nu)$, with $\overline{\mathcal{F}_3(\nu)}=\overline{\nu}$, satisfying  property (Q). 
\end{lemma}

\textbf{Proof. }

We start with a definition.
\vspace{0.5cm}  \begin{definition}
 Let $\nu\in G$. If with  \textbf{y.} relations applied to $\nu$ it is possible to apply one of the relations \textbf{a.}, \textbf{b.}, \textbf{b'.}, \textbf{c.}, \textbf{c'.}, \textbf{d.} or \textbf{x.}, we say $\nu\in L.$\label{d31}
\end{definition}

If with  \textbf{y.} relations applied to $\nu$ it is possible to apply one of the relations \textbf{a.}, \textbf{b.}, \textbf{b'.}, \textbf{c.}, \textbf{c'.}, \textbf{d.} or \textbf{x.}, we do it and we call $\mathcal{F}_3^1(\nu)$ the resulting expression. This is not well defined because there might be many ways of doing this, but we choose one of these ways arbitrarily. We define recursively $\mathcal{F}_3^n(\nu)=\mathcal{F}_3^1(\mathcal{F}_3^{n-1}(\nu)).$

In the following table the symbol $-$ means that the corresponding relation decreases the corresponding $\nu\{-\}$,  the symbol $-0$ means that sometimes it decreases it and sometimes it maintains it equal and the symbol $0$ means it always maintains it equal. By definition $\ndo\in \mathbb{N}^k$, $\ntr, \ncu, \nci\in \mathbb{N}$.

\vspace{1.5cm}
\begin{tabular}{|l|c|c|c|c|}
&$\ndo$&$\ntr $&$\ncu$&$\nci$\\ \hline
 \textbf{c'.}$\ \ \begin{array}{c} s\ \underbracket{s\ r}\\ \underbracket{s\ r}\ s\\r\underbrace{s\ s}\\\hspace{-0.4cm}r\ \ \ s \end{array}= \begin{array}{c} \underbrace{s\ s}r\\\underbracket{s\ \ r}\\r\ \ s \end{array}$ & ---& & &\\ \hline
\textbf{c.}$\ \ \begin{array}{c}s\, \underbracket{r\,\, \, s}\\ \underbrace{s\ s} r \\
\ \  s \hspace{0.4cm}r    \end{array}=\begin{array}{c} \underbracket{s\,\,\, r}\, s\\ \ r\underbrace{s\ s}\\\underbracket{r\ \ s}\\s\ \ r
\end{array}$& ---& & &\\ \hline
\textbf{a.}$\ \ \begin{array}{c} \underbracket{sr}  \\
 \underbracket{rs}\\sr
\end{array} =\begin{array}{c}sr\\ \ \ \\ \ \ 
\end{array}$&---0&---& & \\ \hline 
\textbf{b.} $\ \  \begin{array}{c}
 \underbracket{sr}  \\
 r\hspace{-0.26cm}{_\smallsmile}s\\
\ s
\end{array} =\begin{array}{c}s\,r\hspace{-0.26cm}{_\smallsmile}\\ s \ \,\\ \ \end{array}  $ &---0&---& & \\ \hline
\textbf{b'.}$\ \  \begin{array}{c} \underbracket{r\ s}\\ s\  r\hspace{-0.26cm}{_\smallsmile}\\s\ \ \end{array}= \begin{array}{c}r\hspace{-0.26cm}{_\smallsmile}\ s\\ \ \ s\end{array}$&---0&---& & \\ \hline
\textbf{d.}$\ \ \begin{array}{c}
                                                              s\underbracket{r\ t}\\ \underbracket{s\ t}r\\ t \underbracket{s\ r}\\t\ r\, s
                                                             \end{array}=\begin{array}{c}\underbracket{s\ r}t\\ r\underbracket{s\ t}\\ \underbracket{r\ t}s\\t\ r\, s
                                                             \end{array}$&0&0&---& \\ \hline
\textbf{x.}$\ \ \begin{array}{c} f\\ m\,\, \mathrm{or}\,\,j \end{array}= \begin{array}{c} m\,\, \mathrm{or}\,\,j\\f \end{array}$ &---0&0&---0&--- \\ \hline
\textbf{y.}$\ \ \begin{array}{c} ^if\\ ^jf \end{array}= \begin{array}{c} ^jf\\^if \end{array}$ &0&0&0&0 \\\hline
\end{tabular}
\vspace{1cm}

This table shows that there exists some $N\in \mathbb{N}$ such that $\mathcal{F}_3^N(\nu)\notin L.$ We put $\mathcal{F}_3(\nu)=\mathcal{F}_3^N(\nu)$, and this proves lemma \ref{5.3}
$\hfill \Box$

\subsection{}

\vspace{0.5cm}  \begin{lemma}\label{5.4}
Let $\nu\in G.$ If the $k^{th}$ term of $\mathcal{F}_3(\nu)$ is $^if$ then the $(k+1)^{th}$ term is $^{i+1}f$ or $^{i+1}j.$
\end{lemma}
\textbf{Proof. }
Let $\mathcal{F}_3(\nu)=(^{\omega_k}\hspace{-0.02cm}d_k, \cdots, ^{\omega_1}\hspace{-0.1cm}d_1).$  We use the following notation : if $a<b$ are two natural numbers, then $[a,b]=\{a,a+1,a+2,\ldots, b\}.$\label{d32}

Let $a(1)<b(1)<a(2)<b(2)<\ldots <a(r)<b(r)$ be such that  $$\mathcal{F}_3(\nu)[4]=\bigcup_{i=1}^r \ [a(i),b(i)]$$

Let us fix $1\leq l \leq r$. We only need to prove 
\begin{description}  
\item[A]$d_{b(l)+1}=j, \omega_{b(l)+1}=\omega_{b(l)}+1$ 
\item[B] For $a_l\leq b(l)-m\leq b_l$ we have $\omega_{b(l)-m}=\omega_{b(l)}-m.$
\end{description}
 This means that between the $a(l)^{th}$ term and the $(b(l)+1)^{th}$ term, the picture looks like this : $\begin{array}{c} \underbracket{\ \ \ \  }\\\hspace{0.7cm}\underbracket{\ \ \ \ }\\\hspace{1.9cm}\ddots\\\hspace{3.1cm} \underbracket{\ \ \ \ }\\ \hspace{4.1cm}\underbrace{} \end{array}
$

We start by proving \textbf{A}. Because of the definition of $b(l)$, the $b(l+1)^{th}$ term is a $j$ or $m$ (it can not be an $\alpha$ because $\nu\in G$). Because of the relation \textbf{y.} the $b(l)^{th}$ term does not commute with the $b(l+1)^{th}$ term, so we have four possibilities for the  $b(l+1)^{th}$ term : 

(1) $^{\omega_{b(l)}}m$ \hspace{1.5cm}  (2) $^{\omega_{b(l)}+1}m$   \hspace{1.5cm} (3) $^{\omega_{b(l)}-1}j$   \hspace{1.5cm} (4) $^{\omega_{b(l)}+1}j$

The cases (1), (2) and (3) are impossible respectively by relations \textbf{b.}, \textbf{b'.} and  \textbf{c.}, so we have proved \textbf{A}.

Now we prove \textbf{B} by induction on $m$. We suppose we have proved it for $m$, we will prove it for $m+1.$ We will use the following picture. 

$\begin{array}{c} \hspace{-7cm} a_l\rightarrow\\ \\ \\ \\ \hspace{4.5cm}\cdots \mathbf{x}\ \ \ \mathbf{x}\ \  \blacksquare \ \ \mathbf{a}\  \ \ \ \mathbf{e}\  \ \ \textbf{e}\ \  \cdots \hspace{0.3cm}\mathbf{e}\ \  \ \mathbf{c'}\ \  \mathbf{c}\ \  \ \mathbf{x}\ \ \ \mathbf{x}\ \ \ \mathbf{x}\ \cdots\ \\ \hspace{-3cm}b_l-m\rightarrow \hspace{2cm} \underbracket{\ \ \ \ \ \ \ \ \  }\\\hspace{2cm}\underbracket{\ \ \ \ \ \ \ \ \ }\\\hspace{3.5cm}\ddots\\\hspace{5cm} \underbracket{\ \ \ \ \ \ \ \ \ }\\\hspace{6.3cm} \underbracket{\ \ \ \ \ \ \ \ \ }\\ b_l+1\rightarrow \hspace{6.1cm}\underbrace{\ \ \ \ \ \ } \end{array}
$
In this picture we have drawn from the $(b_l-m)^{th}$ term to the $(b_l+1)^{th}$ term. In the $(b_l-m-1)^{th}$ term we have put some letters, and we will explain what they mean. We put $e=b_l-m-2$. Let $^{\omega_{e}}\hspace{-0.02cm}d_{e}\circ \cdots \circ ^{\omega_1}\hspace{-0.1cm}d_1 :\theta_{s_1}\cdots \theta_{s_n}\rightarrow \theta_{u_1}\cdots \theta_{u_q}$.

 Let us say that a letter $\mathbf{n} \in \{\mathbf{x}, \mathbf{a}, \mathbf{e}, \mathbf{c'}, \mathbf{c}\}$  is in the place that corresponds to $\theta_{u_p}$. Then if the  $(b_l-m-1)^{th}$ term is $^{p-1}f$, by taking this term down  in the picture (using relation \textbf{y} all the times we are allowed), we can use relation \textbf{n}, and by definition of $\mathcal{F}_3(\nu)$, this is not allowed. We conclude that
if  $\blacksquare$  is in the place that corresponds to $\theta_{u_i}$, then the $(b_l-m-1)^{th}$ term of $\mathcal{F}_3(\nu)$ is $^{i-1}f$, and so we prove \textbf{B} and lemma \ref{5.4}.
$\hfill \Box$

Lemma \ref{5.4} joint with the properties of $\mathcal{F}_2$ allows us to finish the proof of proposition \ref{duru}.$\hfill \Box$

\section{Elimination of the $\alpha'$s}
We start this section by introducing some relations that will be useful for proving proposition \ref{alpha}:
\vspace{0.5cm}  \begin{proposition}We have  the following equalities :
\begin{description}
 \item[N1]\label{21} $(\mathrm{id}\otimes f_{r,s})\circ (f_{r,s}\otimes \mathrm{id})\circ
 (\mathrm{id}\otimes \alpha_s)=(\alpha_s\otimes \mathrm{id})  :  \theta_r\rightarrow
 \theta_s\theta_s\theta_r$
\item[N2]\label{22} $  (\mathrm{id}\otimes j_s)\circ (f_{sr}\otimes
 \mathrm{id})=f_{sr}\circ (j_s\otimes \mathrm{id})\circ (\mathrm{id}\otimes f_{r,s}) : \theta_s\theta_r\theta_s
 \rightarrow \theta_r\theta_s $
 \item[N3]\label{23} $(\mathrm{id}\otimes f_{r,s})\circ (f_{r,s}\otimes \mathrm{id})\circ (\mathrm{id}\otimes p_s) =(p_s\otimes \mathrm{id}) \circ f_{r,s}:\theta_r\theta_s\rightarrow
 \theta_s\theta_s\theta_r$
\item[N4]\label{24} $f_{r,s}\circ (\mathrm{id}\otimes \epsilon_s)=\epsilon_s\otimes \mathrm{id} : \theta_r \rightarrow \theta_s\theta_r$

    \item[N5]\label{8}$  (m_r\otimes \mathrm{id})\circ p_r=(\mathrm{id}\otimes m_r)\circ
 p_r=\mathrm{id}   \hspace{2.6cm}   \begin{array}{c} r\\ \overbrace{r\hspace{-0.26cm}{_\smallsmile}\ r}\\\ \ r
 
\end{array}=\begin{array}{c} r\\ \overbrace{r\ r\hspace{-0.26cm}{_\smallsmile}}\\ r\ \ 
 
\end{array}=\begin{array}{c} r\\\ \\ \ 
\end{array}            $
    \item[N6]\label{9}$  j_r\circ (\mathrm{id}\otimes \epsilon_r)=j_r\circ
 (\epsilon_r\otimes \mathrm{id})=\mathrm{id}  \hspace{2.7cm}   \begin{array}{c} r\ \ \\  \underbrace{r\ r\hspace{-0.23cm}{^\smallfrown}}\\ r
 
\end{array}=\begin{array}{c}\ \ r\\  \underbrace{r\hspace{-0.23cm}{^\smallfrown} \ r}\\r
 
\end{array}=  \begin{array}{c} r\\\ \\ \ 
\end{array}    $
    \item[N7]\label{10}$  (\mathrm{id}\otimes j_r)\circ (p_r\otimes \mathrm{id})=$

\vspace{-0.8cm}
$=(j_r\otimes
 \mathrm{id})\circ (\mathrm{id}\otimes p_r)=p_r\circ j_r  \hspace{1cm}       \begin{array}{c} r\ r\\ \overbrace{r\ r} \hspace{-0.5cm} \underbrace{\ \ r}\\ r\ \ \ \ \hspace{-0.2cm}r
 
\end{array}=\begin{array}{c} r\ r\\\underbrace{r\ r}\hspace{-0.5cm}\overbrace{\ \ r}\\\hspace{0.2cm}r\ \ \ \ \hspace{-0.2cm}r
 
\end{array}=\begin{array}{c}\underbrace{r\ r}\\ r\\ \overbrace{r\ r}
 \vspace{0.3cm}
\end{array}       $
\vspace{0.3cm}
\end{description}
\end{proposition}

\textbf{Proof. }
 We will give in order the relations needed to prove each one of this equations (CR means commutation relations) :
\begin{itemize}
 \item \textbf{N1} : \textbf{e}, \textbf{a}
\item \textbf{N2} : $\mathbf{c'}$, \textbf{a}
\item \textbf{N3} : \textbf{e}, \textbf{N2}, CR, \textbf{a}
\item \textbf{N4} : CR, \textbf{e}, $\mathbf{a'}$
\item \textbf{N5} : \ref{7}, \ref{1}
\item \textbf{N6} : CR, \ref{7}, \ref{1}
\item \textbf{N7} : \ref{7}, CR, \ref{7}, \ref{3}, CR
\end{itemize}
$\hfill \Box$

\vspace{0.5cm}  \begin{proposition}\label{alpha}
Let $\tau=\overline{g\circ \xi}$ with $g\in G$, $\xi=^{i}\alpha_s$. There exists a set $\Pi$ and for each $\pi \in \Pi$ a polynomial $\lambda_{\pi}$ and $g_{\pi}\in G$ such that 
$$\tau=\overline{\sum_{\pi\in \Pi}\lambda_{\pi}g_{\pi}} $$
\end{proposition}

\textbf{Proof. }

We start with a lemma

\vspace{0.5cm}  \begin{lemma}\label{proz} For proving  proposition \ref{alpha} it suffices to prove it for  $\xi=^{i}p_s$ and $\xi=^{i}\epsilon_s$. 
\end{lemma}
\textbf{Proof. }
We have that $\mu_1=(\mathcal{F}_3(g), ^{i}\alpha_s)$ is an expression of the  morphism
 $\tau$. With commutation relations we change $\mu_1$ in $\mu_2$, an
 expression where the $\alpha$ is as far as possible to the left (in the picture is in the lowest place possible). Let us say that in
 $\mu_2$ we have that $^{p}\alpha_s$  is the $k^{th}$ term. We have seven possibilities for the $(k+1)^{th}$ term : $(1)\,\, ^{p-1}j $, $(2)\,\,
 ^{p}j $, $(3)\,\, ^{p+1}j$, $(4)\,\, ^{p}m $, $(5)\,\, ^{p+1}m $,
 $(6)\,\, ^{p+1}f $ and $(7)\,\, ^{p-1}f.$

\begin{itemize}\vspace{0.3cm}
\item In the case (1), by  relation \ref{7} we are in the case
 $\xi=^{p}p_s$. 
\vspace{0.3cm}\item In the case  (2), by relation \ref{5} we find that 
 $\tau=0$. 
\vspace{0.3cm}\item In the case  (3) by definition of $p_s$  we are in the case $\xi=^{p}\hspace{-0.13cm}p_s$. 
\vspace{0.3cm}\item In the case (4)  by relation \ref{6} we are in the case $\xi=^{p}\hspace{-0.13cm}\epsilon_s$. 
\vspace{0.3cm}\item In the case (5) by definition of $\epsilon_s$ we are in the case $\xi=^{p}\hspace{-0.13cm}\epsilon_s$.
\vspace{0.3cm}\item In the case (6), by remark \ref{ya} we have that there exists $c$
such that the $(k+r)^{th}$ term is  $^{p+r}f$ for every $0\leq r< c$ and 
 $^{p+c}j$ for $r=c.$ But with the relation \textbf{N2} we can easily verify the following equality by induction  : $$^{p+c}j\circ
 ^{p+c-1}f\circ \ldots \circ ^{p+1}f=(^{p+c-1}f\circ ^{p+c-2}f\circ \cdots
 \circ ^{p+1}f)\circ ^{p+1}j\circ (^{p+2}f\circ ^{p+2}f\circ
 \cdots \circ ^{p+c}f) $$
So if we replace the left hand side of this equation by it is right hand side we can take the  $\alpha$ more to the left, and we find a new expression of  $\tau$ : $(g_1, ^{p+1}\hspace{-0.15cm}j,
 ^p\hspace{-0.05cm}\alpha, g_2)$, with $g_1,g_2\in G$ and by the definition of $p_s$ we get to the case $\xi=^{p}p_s$.
\vspace{0.3cm}\item In the case (7), by the remark \ref{ya} we have that the
 $(k+2)^{th}$ term is $^pf$. So if we apply  relation
 \textbf{N1} we arrive to a new expression $\mu_3$ of $\tau$
 : $\mu_3=(g_1, ^{p-1}\hspace{-0.1cm}\alpha_s, g_2)$, $g_1, g_2\in G$, where $g_1$ is a good $g$-expression, the $\alpha$ is still in the $k^{th}$ term and $\mu_3$ has strictly less terms than $\mu_2.$
    Now we repeat the process of taking in $\mu_3$ the 
 $\alpha$ as far to the left as possible and if we arrive another time to the case  (7), we find a corresponding $\mu_4$. If we repeat this process enough times, finally we will arrive to one of the other 6 cases. This finishes the proof of the lemma.
\end{itemize}
$\hfill \Box$

\textbf{Proof of the proposition \ref{alpha} for $\xi=^{i}p_s$}. The proof of this case is very similar to the proof of lemma \ref{proz}, but we use different relations. We have that $\mu_1=(\mathcal{F}_3(g),^{i}\hspace{-0.17cm}p_s)$ is an expression of the morphism
 $\tau$. With commutation relations we change $\mu_1$ in $\mu_2$, an expression where the  $p$ is as far to the left as possible. Let us say that in 
 $\mu_2$ we have that the $k^{th}$ term is  $^{p}p_s$. We have seven possibilities for the  $(k+1)^{th}$ term : $(1)\,\, ^{p-1}j $, $(2)\,\,
 ^{p}j $, $(3)\,\, ^{p+1}j$, $(4)\,\, ^{p}m $, $(5)\,\, ^{p+1}m $,
 $(6)\,\, ^{p+1}f $ and $(7)\,\, ^{p-1}f.$
\begin{itemize}
\vspace{0.3cm}\item In the cases  (1) and (3) we use  relation \textbf{N7}, and we have that the $p$ is more to the left than before.
\vspace{0.3cm}\item In the case (2) we have that  $\tau=0$ because the relations \ref{3} and \ref{5} tells us that  $j_s\circ p_s=0.$ 
\vspace{0.3cm}\item In the cases (4) and (5) we use  relation \textbf{N5} and we find an expression of $\tau$ that is in $G$.
\vspace{0.3cm}\item In the case (6), by a similar argument to that of case (6) of  lemma \ref{proz}, we go back to case (3).
\vspace{0.3cm}\item In the case (7) by a similar argument to that of case  (7) of lemma \ref{proz},  but using the relation  \textbf{N3} instead of relation  \textbf{N1} we see that, as in the cases (1), (3) and (6), the $p$ is more to the left than before. So if we repeat enough times we will go back to one of the cases that are left, this means, cases (2), (4) or (5).
\end{itemize}

\vspace{0.5cm}
\textbf{Proof of the proposition \ref{alpha} for $\xi=^{i}\epsilon_s$}. We have that $\mu_1=\mathcal{F}_3(g)\circ ^{i}\epsilon_s$ is an expression of  morphism
 $\tau$. With commutation relations we change $\mu_1$ in $\mu_2$, an 
 expression where the $\epsilon$ is as far to the left as possible. Let us say that in $\mu_2$ we have that the $k^{th}$ term is $^{p}\epsilon_s$. We have five possibilities for the  $(k+1)^{th}$ term : $(1)\,\, ^{p-1}j $, $(2)\,\,
 ^{p}j $, $(3)\,\, ^{p}m $,
 $(4)\,\, ^{p}f $ and $(5)\,\, ^{p-1}f.$
\begin{itemize}
\vspace{0.3cm}\item In the cases (1) and (2), using relation \textbf{N6} we find an expression of $\tau$ that is in $G$.
\vspace{0.3cm}\item In the case (4), by a similar argument to that of the case (6) in lemma \ref{proz}, we go back to case (2).
\vspace{0.3cm}\item In the case (5), by a similar argument to that of the case (7) in lemma \ref{proz}, but using relation  \textbf{N4}  instead of relation \textbf{N1} we see that the $\epsilon$ is more to the left than before. So, if we repeat enough times we will go back to one of the cases that are left.
\vspace{0.3cm}\item  Case (3) is treated  in section \ref{tak}.
\end{itemize}
$\hfill \Box$

By using proposition \ref{alpha} repeatedly we have the following

\vspace{0.5cm}  \begin{corollary}\label{seis}
 Let $\nu$ be an $R$-expression. There exists a set $\Pi$ and for each $\pi \in \Pi$ a polynomial $\lambda_{\pi}$ and $g_{\pi}\in G$ such that 
$$\overline{\nu}=\overline{\sum_{\pi\in \Pi}\lambda_{\pi}g_{\pi}} $$
\end{corollary}

\section{Good order}
The purpose of this section is to change a good $g$-expression into a good $g$-expression in good order. 
\vspace{0.5cm}  \begin{proposition}
Let $t>t'$ and $\alpha, \beta\in \{m, ch, cch\}$. There exists $u$ and $u'$ such that $\alpha(t)\circ \beta(t')=\beta(u')\circ \alpha(u)$.
\end{proposition}
\textbf{Proof. }

$\bullet$  Let us consider the case $\beta=\mathrm{ch}$ or $\mathrm{cch}$ and $\alpha=m$. We have two possible cases. In the first one $m(t)$ commutes with $\beta(t')$, and in the second one they do not commute, so we use commutation relations and relation \textbf{b.} 
\begin{itemize}
 \item  $ m(t)\circ \beta(t')=\beta(t')\circ m(t+1)$

$\begin{array}{c} \underbracket{\ \ \ \  }\\\hspace{0.7cm}\underbracket{\ \ \ \ }\\\hspace{1.9cm}\ddots\\\hspace{3.1cm} \underbracket{\ \ \ \ }\\ \hspace{4.1cm}\underbrace{}\\  \hspace{-2.2cm}{_\smallsmile}\end{array} \rightarrow   \begin{array}{c}\hspace{-2.2cm}{_\smallsmile}\\ \underbracket{\ \ \ \  }\\\hspace{0.7cm}\underbracket{\ \ \ \ }\\\hspace{1.9cm}\ddots\\\hspace{3.1cm} \underbracket{\ \ \ \ }\\ \hspace{4.1cm}\underbrace{}\end{array}$  
\vspace{0.8cm}
\item   $ m(t)\circ \beta(t')=\beta(t')\circ m(t)$

$\begin{array}{c} \underbracket{\ \ \ \  }\\\hspace{0.7cm}\underbracket{\ \ \ \ }\\\hspace{1.9cm}\ddots\\\hspace{3.1cm} \underbracket{\ \ \ \ }\\ \hspace{4.1cm}\underbrace{}\\  \hspace{1cm}{_\smallsmile}\end{array} \rightarrow   \begin{array}{c}\hspace{1.8cm}{_\smallsmile}\\ \underbracket{\ \ \ \  }\\\hspace{0.7cm}\underbracket{\ \ \ \ }\\\hspace{1.9cm}\ddots\\\hspace{3.1cm} \underbracket{\ \ \ \ }\\ \hspace{4.1cm}\underbrace{}\end{array}$ 

\end{itemize}

\vspace{0.5cm}
$\bullet$ The case $i=1, 2$ or $3$ and $k=1$ is easy because $\alpha(t)$ commutes with $m(t').$

\vspace{0.5cm}
$\bullet$ The last case is $i=2$ or $3$ and $k=2$ or $3$. We will only treat the case $i=2, k=2,$ the other ones are similar.
We will prove that $\mathrm{ch}(t)\circ \mathrm{ch}(t')=\mathrm{ch}(t')\circ \mathrm{ch}(t).$ For this we need two preliminary lemmas.

\vspace{0.5cm}  \begin{lemma}\label{7.2}
 If $f=^{i_q}f\circ \cdots \circ ^{i_1}f,$\  $g=^{k_p}f\circ \cdots \circ ^{i_1}f\in \mathrm{Hom}(\theta_{s_1}\cdots \theta_{s_n},\theta_{u_1}\cdots \theta_{u_l})$, then $f=g.$
\end{lemma}
\textbf{Proof. }
Relations \textbf{a.} and \textbf{d.} are exactly the relations defining the symmetric group.$\hfill \Box$

\vspace{0.5cm}  \begin{lemma}\label{7.3}
 The following expressions represent the same morphism : $(j^i, f^{i-2}, f^{i-1}, f^i)$ and $(f^i, f^{i+1}, j^{i+1}).$
\end{lemma}
\textbf{Proof. }
In the following chain of relations we apply respectively commutation relations, relation \textbf{c.} and relation \textbf{a.} :

$\begin{array}{c} \underbracket{\ \ \ \ \ } \\ \hspace{0.8cm}\underbracket{\ \ \ \ \ } \\\hspace{1.7cm}\underbracket{\ \ \ \ \ } \\\underbrace{}  \end{array} \rightarrow
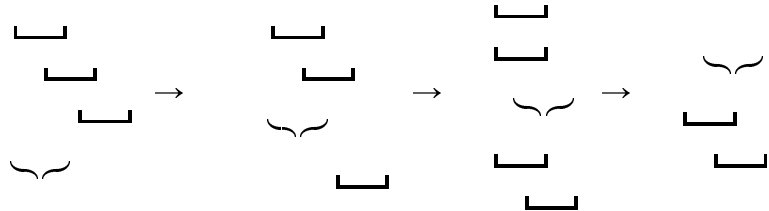
\begin{array}{c} \underbracket{\ \ \ \ \ } \\\hspace{0.8cm}\underbracket{\ \ \ \ \ } \\ \underbrace{}\\ \hspace{1.7cm}\underbracket{\ \ \ \ \ } \end{array} \rightarrow
\begin{array}{c}\underbracket{\ \ \ \ \ } \\\underbracket{\ \ \ \ \ } \\ \hspace{0.6cm}\underbrace{}\\ \underbracket{\ \ \ \ \ } \\\hspace{0.8cm}\underbracket{\ \ \ \ \ }  \end{array} \rightarrow
\begin{array}{c} \hspace{0.6cm}\underbrace{}\\ \underbracket{\ \ \ \ \ } \\\hspace{0.8cm}\underbracket{\ \ \ \ \ } \end{array}$

$\hfill \Box$

Now we are able to finish the case $i=k=2$. In the following chain of equalities we apply respectively commutation relations to take an $j$ down, then  lemma \ref{7.2} for reordering a composition of $f'$s,  and finally with commutation relations we take up an $j$ and apply lemma \ref{7.3} :

$\begin{array}{c} \underbracket{\ \ \ \  \ }\\\hspace{0.8cm}\underbracket{\ \ \ \ \ }\\\hspace{2.2cm}\ddots\\\hspace{3.5cm} \underbracket{\ \ \ \ \ }\\ \hspace{4.5cm}\underbrace{} \\ \hspace{-3cm} \underbracket{\ \ \ \ \  }\\\hspace{-2.2cm}\underbracket{\ \ \ \ \ }\\\hspace{-0.8cm}\ddots\\\hspace{0.5cm} \underbracket{\ \ \ \ \ }\\ \hspace{1.5cm}\underbrace{}  \end{array}\rightarrow
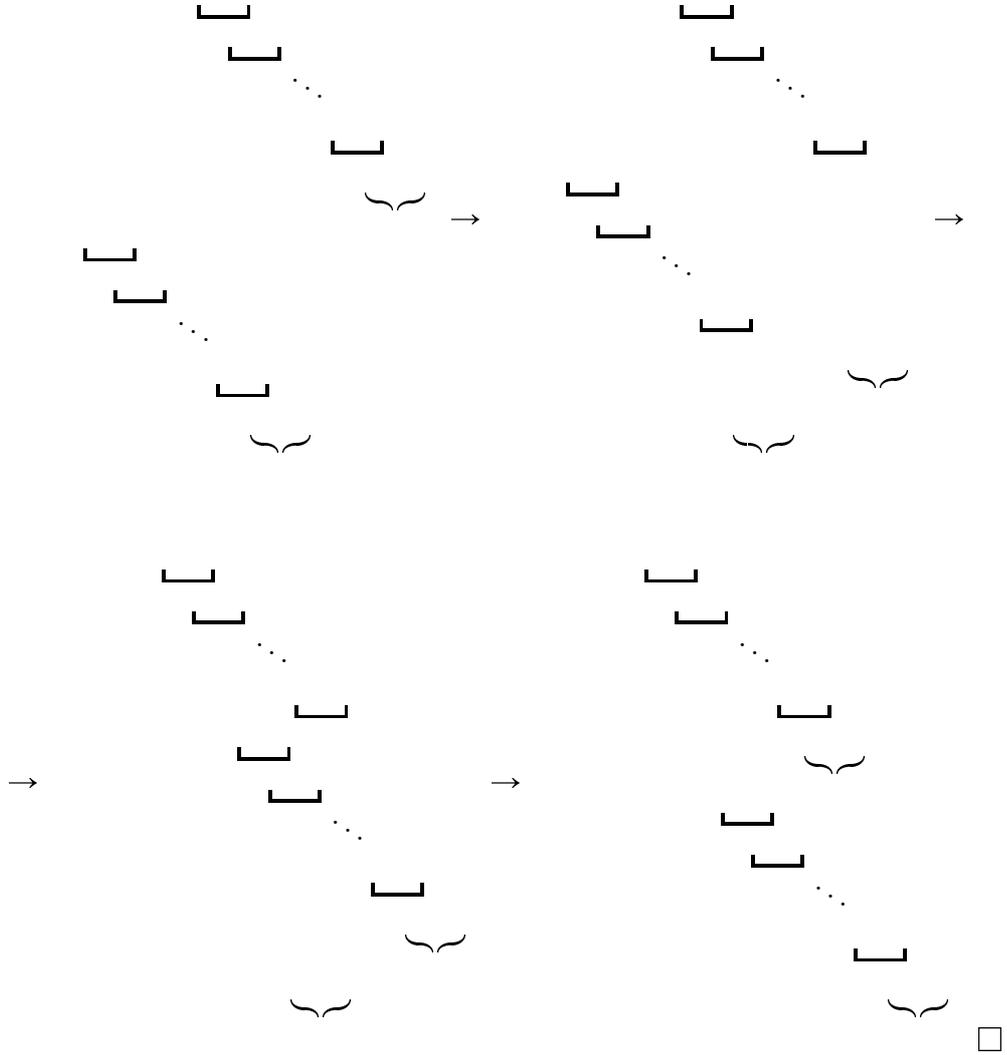
\begin{array}{c}  \underbracket{\ \ \ \  \ }\\\hspace{0.8cm}\underbracket{\ \ \ \ \ }\\\hspace{2.2cm}\ddots\\\hspace{3.5cm} \underbracket{\ \ \ \ \ } \\ \hspace{-3cm} \underbracket{\ \ \ \ \  }\\\hspace{-2.2cm}\underbracket{\ \ \ \ \ }\\\hspace{-0.8cm}\ddots\\\hspace{0.5cm} \underbracket{\ \ \ \ \ }\\ \hspace{4.5cm}\underbrace{}\\ \hspace{1.5cm}\underbrace{}  \end{array}\rightarrow$

\vspace{1cm}
$\rightarrow
\begin{array}{c} \hspace{-2cm} \underbracket{\ \ \ \ \  }\\\hspace{-1.2cm}\underbracket{\ \ \ \ \ }\\\hspace{0.2cm}\ddots\\\hspace{1.5cm} \underbracket{\ \ \ \ \ }\\ \underbracket{\ \ \ \  \ }\\\hspace{0.8cm}\underbracket{\ \ \ \ \ }\\\hspace{2.2cm}\ddots\\\hspace{3.5cm} \underbracket{\ \ \ \ \ } \\ \hspace{4.5cm}\underbrace{}\\ \hspace{1.5cm}\underbrace{}\end{array}\rightarrow
\begin{array}{c}  \hspace{-2cm} \underbracket{\ \ \ \ \  }\\\hspace{-1.2cm}\underbracket{\ \ \ \ \ }\\\hspace{0.2cm}\ddots\\\hspace{1.5cm} \underbracket{\ \ \ \ \ }\\ \hspace{2.3cm}\underbrace{}\\ \underbracket{\ \ \ \  \ }\\\hspace{0.8cm}\underbracket{\ \ \ \ \ }\\\hspace{2.2cm}\ddots\\\hspace{3.5cm} \underbracket{\ \ \ \ \ } \\ \hspace{4.5cm}\underbrace{} \end{array} $

$\hfill \Box$

\vspace{0.5cm}  \begin{definition}
 If $\nu$ is a good $g$-expression we will name $\mathcal{F}_4(\nu)$\label{d29} the good $g$-expression in the good order such that $\overline{\mathcal{F}_4(\nu)}=\overline{\nu}$
\end{definition}

\section{Conclusion of the proof of theorem \ref{bacanz}}\label{section8}

\subsection{}
\vspace{0.5cm}  \begin{definition}
 If $\tau=\sum_i g_i\lambda_i$, with $I$ a finite set, $\lambda_i$ polynomials and $g_i\in G $, we define $\mathcal{F}_3(\tau)=\sum_i \mathcal{F}_3(g_i)\lambda_i$ and $\mathcal{F}_4(\tau)=\sum_i \mathcal{F}_4(g_i)\lambda_i$.
\end{definition}

\vspace{0.5cm}  \begin{definition}
 Let $\nu\in G.$ We define $\mathcal{F}_5^1(\nu)=\mathcal{F}_4\mathcal{F}_3\mathcal{F}_2(\nu).$  This is well defined, because $\mathcal{F}_2(\nu)$ is a linear combination of elements of $G$ and $\mathcal{F}_3(\nu')$ is a good $g$-expression. 

We define inductively  $\mathcal{F}^n_5(\nu)$ :
If $\mathcal{F}_5^n(\nu)\ce\neq 0$ we define $\mathcal{F}_5^{n+1}(\nu)=\mathcal{F}_5^1\mathcal{F}_5^n(\nu)$ and if $\mathcal{F}_5^n(\nu)\ce=0,$ we define $\mathcal{F}_5^{n+1}(\nu)=\mathcal{F}_5^n(\nu).$ 
\end{definition}

\vspace{0.5cm}  \begin{proposition}\label{socho}\label{d30}
 $\mathcal{F}_5^n(\nu)$ stabilizes for $n$ large.
\end{proposition}
\textbf{Proof. }
 Let us suppose it doesn't stabilize. This means that for all $n\in \mathbb{N},$ we have $\mathcal{F}_5^n(\nu)\ce\neq 0.$ So we apply infinitely many times relation \ref{12} in this process. When we apply relation \ref{12} to an expression we obtain a sum of three expressions, the first one  is $\begin{array}{c} r\hspace{-0.26cm}{_\smallsmile}\ r\\ \ \ r 
\end{array}$, the second one $\begin{array}{c}\underbrace{r\ r}\\ r\\ \ r^{\diamond}
\end{array}$, and the third one is $\begin{array}{c}\underbrace{r\ r}\\ r\\  ^{\diamond}r\ \  
\end{array}$. In this three expressions  relation \ref{12} always decreases strictly $\non\in \mathbb{N}$ and the other relations used in $\mathcal{F}_{2},$ $\mathcal{F}_{3}$ and $\mathcal{F}_{4}$ don't change $\non$. To see this, we make a brief review of all relations used in defining 
$\mathcal{F}_{2},$ $\mathcal{F}_{3}$ and $\mathcal{F}_{4}$ : 
\begin{itemize}
 \item $\mathcal{F}_{2} : $ commutation relations, \textbf{a.} in the opposite sense, \ref{12}, \textbf{g}, \textbf{f} 
\item $\mathcal{F}_{3} : $ \textbf{a.}, \textbf{b.}, \textbf{b'.}, \textbf{c.}, \textbf{c'.}, \textbf{d.}, \textbf{x.}, \textbf{y.}
\item $\mathcal{F}_{4} : $ commutation relations,   \textbf{a.},\textbf{b.}, \textbf{c.}, \textbf{d.}
\end{itemize}
 So we have a contradiction, which allows us to conclude the proof.
$\hfill \Box$

\vspace{0.5cm}  \begin{remark}\label{seche}
 As $\mathcal{F}_4(\nu)$ is a good $g$-expression in the good order, $\mathcal{F}_5(\nu)$ is a good $g$-expression in the good order satisfying $\mathcal{F}_5(\nu)\ce=0.$
\end{remark}

\subsection{}

We start by a useful reformulation of proposition \ref{FL} in our terminology :

\vspace{0.5cm}  \begin{corollary}
 Let $(W,\mathcal{S})$ be a right-angled Coxeter system. Let $(s_1,\ldots,s_n)\in\mathcal{S}^n.$ We have the equality of sets :

$FL(s_1,\ldots, s_n)=\{\nu$ good $g$-expressions in $\mathrm{Hom}(\theta_{s_1}\cdots \theta_{s_n},R)$ in good order,\newline 
\hspace*{6cm} with $\nce=\nun=0\}$
\end{corollary}

In section \ref{rest} we showed that in order to prove Theorem \ref{bacanz} we only need to prove the following proposition :

\vspace{0.5cm}  \begin{proposition}
 Let $(s_1,\ldots,s_n)\in \mathcal{S}^n.$ If $f\in \mathrm{Hom}(\theta_{s_1}\cdots \theta_{s_n},R)$ then there exists $I$ a finite set and for $i\in I$, $a_i\in FL(s_1,\ldots,s_n)$ and $\lambda_i$ polynomials such that $f=\sum_ia_i\lambda_i.$
\end{proposition}
As $f$ is a linear combination over $A$ of $R$-expressions, we can restrict our attention to the case where $f$ is an $R$-expression. By corollary \ref{seis} we can restrict to the case $f\in G.$ By proposition \ref{socho} and remark \ref{seche} we can restrict to the case where $f$ is a good $g$-expression in the good order and with $\nce=0.$ So, the following lemma allows us to conclude the proof of Theorem \ref{bacanz} :

\vspace{0.5cm}  \begin{lemma}
 If $\nu$ is a good $R$-expression in the good order and $\nce=0$, then $\nun=0.$
\end{lemma}
\textbf{Proof. }
Suppose $\nun\neq 0.$ Let us put  $k=\ndi$. Let $^zj_s$ be the $k^{\mathrm{th}}$ term of $\nu$. The $(k+1)^{\mathrm{th}}$ term of $\nu$ can not be $^zm_s$, because $\nce=0$.

Let us define by induction the natural number $N_p$ for $r\geq p\geq 1$ (where $r$ will be defined in the process). We define $N_1=z+1$. Let us suppose that we have defined $N_{p}$.  \vspace{0.3cm}
\begin{itemize}
\vspace{0.3cm}\item If the $(k+p)^{\mathrm{th}}$ term of $\nu$ is $^im$, then $N_{p+1}=N_p$. 
\vspace{0.3cm}\item Suppose that the $(k+p)^{\mathrm{th}}$ term of $\nu$ is $^ij$. If $i>N_p-1$ then $N_{p+1}=N_p$ and if $i=N_p-1$ then $r=p$.  
\vspace{0.3cm}\item Suppose that the $(k+p)^{\mathrm{th}}$ term of $\nu$ is $^if$. If $i\notin \{N_p-1, N_p-2\}$ then $N_{p+1}=N_p$. If $i=N_{p-1}$ then $N_{p+1}=N_p+1$ and if $i=N_{p-2}$ then $N_{p+1}=N_p-1$.
\end{itemize}

The fact that $\nu$ is a good $R$-expression in the good order allows us to conclude that $N_p$ is well defined.
By construction, the $N_p^{\mathrm{th}}$ element of  $\widehat{\nu}_{k+p}$ is always the same element of $\mathcal{S}$, it is  of left type, and it is evident that $r$ is a finite number.
This contradicts the definition of $\ndi$, so we have a contradiction. We conclude that $\nun= 0.$ 
$\hfill \Box$



\end{document}